\documentclass[a4paper,fleqn]{cas-dc}
\usepackage[numbers]{natbib}
\usepackage{times}
\usepackage{hyperref}

\usepackage{amsthm}
\usepackage{amssymb}
\usepackage{amsmath}
\usepackage{algorithm}
\usepackage{algorithmic}
\usepackage{graphicx}
\usepackage{caption}
\newtheorem{assumption}{Assumption}
\newtheorem{definition}{Definition}
\newtheorem{problem}{Problem}
\newtheorem{lemma}{Lemma}
\newtheorem{theorem}{Theorem}
\newtheorem*{proof*}{Proof}
\usepackage[utf8]{inputenc}
\usepackage[numbers]{natbib}
\usepackage{url}
\usepackage{amssymb}
\usepackage{amsmath}
\usepackage{graphicx}
\newcommand{\tabincell}[2]{\begin{tabular}{@{}#1@{}}#2\end{tabular}}
\usepackage{multirow}
\usepackage{algorithm}
\usepackage{algorithmic}
\usepackage{setspace}

\begin{document}
\shortauthors{Junxiang Wang, Hongyi Li, and Liang Zhao.}
\shorttitle{Accelerated Gradient-free Neural Network Training by Multi-convex Alternating Optimization}
\title[mode = title]{Accelerated Gradient-free Neural Network Training by Multi-convex Alternating Optimization} 
\author[1]{Junxiang Wang}
\ead{jwan936@emory.edu}
\author[2]{Hongyi Li}
\ead{lihongyi@stu.xidian.edu.cn}
\author[1]{Liang Zhao}
\ead{lzhao41@emory.edu}

\address[1]{Emory University, 201 Dowman Dr, Atlanta, GA, USA, 30322}
\address[2]{The State Key Laboratory of Integrated Service Network, Xidian University, Xi'an, Shannxi, China, 710071}

\begin{abstract}
In recent years, even though Stochastic Gradient Descent (SGD) and its variants are well-known for training neural networks, it suffers from limitations such as the lack of theoretical guarantees, vanishing gradients, and excessive sensitivity to input. To overcome these drawbacks, alternating minimization methods have attracted fast-increasing attention recently. As an emerging and open domain, however, several new challenges need to be addressed, including 1) Convergence properties are sensitive to penalty parameters, and 2) Slow theoretical convergence rate. We, therefore, propose a novel monotonous Deep Learning Alternating Minimization (mDLAM) algorithm to deal with these two challenges. Our innovative inequality-constrained formulation infinitely approximates the original problem with non-convex equality constraints, enabling our convergence proof of the proposed mDLAM algorithm regardless of the choice of hyperparameters. Our mDLAM algorithm is shown to achieve a fast linear convergence by the Nesterov acceleration technique. Extensive experiments on multiple benchmark datasets demonstrate the convergence, effectiveness, and efficiency of the proposed mDLAM algorithm. Our code is available at \url{https://github.com/xianggebenben/mDLAM}.
\end{abstract}

\maketitle
\section{Introduction}
\indent Stochastic Gradient Descent (SGD) and its variants have become popular optimization methods for training deep neural networks. Many variants of SGD methods have been presented, including  SGD momentum \cite{sutskever2013importance},  AdaGrad \cite{duchi2011adaptive}, RMSProp \cite{tielemandivide}, Adam \cite{kingma2014adam} and AMSGrad \cite{j.2018on}.While many researchers have provided solid theoretical guarantees on the convergence of SGD \cite{kingma2014adam,j.2018on,sutskever2013importance}, the assumptions of their proofs cannot be applied to problems involving deep neural networks, which are highly nonsmooth and nonconvex.  Aside from the lack of theoretical guarantees, several additional drawbacks restrict the applications of SGD. It suffers from the gradient vanishing problem, meaning that the error signal diminishes as the gradient is backpropagated, which prevents the neural networks from utilizing further training \cite{taylor2016training}, and the gradient of the activation function is highly sensitive to the input (i.e. poor conditioning), so a small change in the input can lead to a dramatic change in the gradient.\\
\indent To tackle these intrinsic drawbacks of gradient descent optimization methods, alternating minimization methods have started to attract attention as a potential way to solve deep learning problems. A neural network problem is reformulated as a nested function associated with multiple linear and nonlinear transformations across multi-layers. This nested structure is then decomposed into a series of linear and nonlinear equality constraints by introducing auxiliary variables and penalty hyperparameters. The linear and nonlinear equality constraints generate multiple subproblems, which can be minimized alternately. Many recent alternating minimization methods have focused on applying the Alternating Direction Method of Multipliers (ADMM) \cite{taylor2016training,wang2019admm}, Block Coordinate Descent (BCD) \cite{pmlr-v97-zeng19a} and Method of Auxiliary Coordinates (MAC) \cite{carreira2014distributed} to replace a nested neural network with a constrained problem without nesting, with empirical evaluations demonstrating good scalability in terms of the number of layers and high accuracy on the test sets. These methods also avoid gradient vanishing problems and allow for non-differentiable activation functions such as binarized neural networks \cite{courbariaux2015binaryconnect}, as well as allowing for complex non-smooth regularization and the constraints that are increasingly important for deep neural architectures that are required to satisfy practical requirements such as interpretability, energy-efficiency, and cost awareness \cite{carreira2014distributed}. The ADMM, as a representative of alternating minimization methods, has been explored extensively for different neural network architectures. It was first used to solve the Multi-Layer Perceptron (MLP) problem with convergence guarantees \cite{taylor2016training,wang2019admm}, and then was extended to other architectures such as Recurrent Neural Network (RNN) \cite{yu2021admmirnn}. Recently, it was utilized to achieve parallel neural network training \cite{wang2020toward,wang2021quantized}.
\\
\indent However, as an emerging domain, alternating minimization for deep model optimization suffers from several unsolved challenges including \textbf{1.  Convergence properties are sensitive to penalty parameters.}
One recent work by Wang et al. firstly proved the convergence guarantee of ADMM in the MLP problem \cite{wang2019admm}. However, such convergence guarantee is dependent on the choice of penalty hyperparameters:  the convergence cannot be guaranteed anymore when penalty hyperparameters are small; \textbf{2. Slow convergence rate.} To the best of our knowledge, almost all existing alternating minimization methods can only achieve a sublinear convergence rate. For example, the convergence rate of the ADMM and the BCD is proven to be $O(1/k)$, where $k$ is the number of iterations \cite{wang2019admm,wang2022multi,pmlr-v97-zeng19a}. Therefore, there is still a lack of a theoretical framework that can achieve a faster convergence rate.\\
 \indent To simultaneously address these technical problems,  we propose a new formulation of the neural network problem, along with a novel monotonous Deep Learning Alternating Minimization (mDLAM) algorithm. Specifically, we, for the first time, transform the original neural network optimization problem into an inequality-constrained problem that can infinitely approximate the original one. Applying this innovation to an inequality-constraint-based transformation not only ensures the convexity and hence easily ensures global minima of all subproblems, but also prevents the output of a nonlinear function from changing much and reduces sensitivity to the input. Moreover,  our proposed mDLAM algorithm can achieve a linear convergence rate theoretically, and the choice of hyperparameters does not affect the convergence of our mDLAM algorithm theoretically. Extensive experiments on four benchmark datasets show the convergence, effectiveness, and efficiency of the proposed mDLAM algorithm. 
Our contributions in this paper include:
\begin{itemize}
\item  We propose a novel formulation for neural network optimization. The deeply nested activation functions are disentangled into separate functions innovatively coordinated by inherently convex inequality constraints.

\item We present an efficient optimization algorithm. A quadratic approximation technique is utilized to avoid matrix inversion. Every subproblem has a closed-form solution. The Nesterov acceleration technique is applied to further boost convergence.
\item  We investigate the convergence of the proposed mDLAM algorithm under mild conditions. The new mDLAM algorithm is guaranteed to converge to a stationary point whatever hyperparameters we choose. Furthermore, the proposed mDLAM algorithm is shown to achieve a linear convergence rate, which is faster than existing methods.
\item Extensive experiments have been conducted to demonstrate the effectiveness of the proposed mDLAM algorithm. We test our proposed mDLAM algorithm on four benchmark datasets. Experimental results illustrate that our proposed mDLAM algorithm is linearly convergent on four datasets, and outperforms consistently state-of-the-art optimizers. Sensitivity analysis on the running time shows that it increases linearly with the increase of hidden units and hyperparameters.
\end{itemize}
\indent  The rest of this paper is organized as follows: In Section \ref{sec:related work}, we summarize recent related research work to this paper. In Section \ref{sec:algorithm}, we formulate the MLP training problem and present the proposed mDLAM algorithm to train the MLP model. Section \ref{sec:convergence} details  convergence properties of the proposed mDLAM algorithm. Extensive experiments on benchmark datasets are shown in Section \ref{sec:experiment}, and Section \ref{sec:conclusion} concludes this work.
\section{Related Work}
\label{sec:related work}
\indent  All existing works on deep learning optimization methods fall into two major categories: SGD methods and alternating minimization methods, which are shown as follows:\\
\indent \textbf{SGD methods:} The renaissance of SGD can be traced back to 1951 when Robbins and Monro published the first paper \cite{robbins1951textordfemininea}. The famous back-propagation algorithm was introduced by Rumelhart et al. \cite{rumelhart1986learning}. Many variants of SGD methods have since been presented, including the use of Polyak momentum, which accelerates the convergence of iterative methods \cite{polyak1964some}, and research by Sutskever et al., who highlighted the importance of Nesterov momentum and initialization \cite{sutskever2013importance}. During the last decade, many well-known SGD methods which are incorporated with adaptive learning rates have been proposed by the deep learning community, which include but are not limited to AdaGrad \cite{duchi2011adaptive}, RMSProp \cite{tielemandivide}, Adam \cite{kingma2014adam}, AMSGrad \cite{j.2018on}, Adabelief \cite{zhuang2020adabelief} and Adabound \cite{luo2018adaptive}.\\
\indent \textbf{Applications of alternating minimization methods for deep learning:} Many recent works have applied alternating minimization algorithms to specific deep learning applications. For example, Taylor et al. and Wang et al. presented the ADMM to solve an MLP training problem via transforming it into an equality-constrained problem, where many subproblems split by ADMM can be solved efficiently \cite{taylor2016training,wang2019admm}, Wang et al. proposed a parallel ADMM algorithm to train deep MLP models \cite{wang2020toward}, and a similar algorithm was extended to Graph Augmented-MLP (GA-MLP) models with the introduction of the quantization technique \cite{wang2021quantized}. Zhang et al. handled Very Deep Supervised Hashing (VDSH) problems by utilizing an ADMM algorithm to overcome issues related to vanishing gradients and poor computational efficiency \cite{zhang2016efficient}. Zhang and Bastiaan trained a deep neural network by utilizing ADMM with a graph \cite{zhang2017training} and Askari et al. introduced a new framework for MLP models and optimize the objective using BCD methods \cite{askari2018lifted}. Li et al. proposed an ADMM algorithm to achieve distributed learning of Graph Convolutional Network (GCN) via community detection \cite{li2021community}. Qiao et al. proposed an inertial proximal alternating minimization to train MLP models \cite{qiao2021inertial}.\\
\indent \textbf{Convergence of alternating minimization methods for deep learning:}
Aside from applications, the other branch of works mathematically proves the convergence of the proposed alternating minimization approaches. For instance,  Carreira and Wang proposed a method involving the use of auxiliary coordinates to replace a nested neural network with a constrained problem without nesting \cite{carreira2014distributed}. Lau et al. proposed a BCD optimization framework and proved the convergence via the Kurdyka-Lojasiewicz (KL) property \cite{lau2018proximal}, while Choromanska et al. proposed a BCD algorithm for training deep MLP models based on the concept of co-activation memory \cite{choromanska2018beyond}, and a BCD algorithm with R-linear convergence was proposed by Zhang and Brand to train Tikhonov regularized deep neural networks \cite{zhang2017convergent}. Jagatap and Hegde introduced a new family of alternating minimization methods and prove their convergence to a global minimum \cite{jagatap2018learning}. Yu et al. proved the convergence of the proposed ADMM for RNN models \cite{yu2021admmirnn}. However, to the best of our knowledge, there is a lack of a flexible framework which allows for different activation functions and guarantees a linear convergence rate.

\section{Model and Algorithms}
\label{sec:algorithm}
\begin{table}
\scriptsize
 \centering
 \begin{tabular}{cc}
 \hline
 Notations&Descriptions\\ \hline
 $L$& Number of layers.\\
 $W_l$& The weight vector in the $l$-th layer.\\
 $z_l$& The output of the linear mapping in the $l$-th layer.\\
 $h_l(z_l)$& The nonlinear activation function in the $l$-th layer.\\
 $a_l$& The output of the $l$-th layer.\\
 $x$& The input matrix of the neural network.\\
 $y$& The predefined label vector.\\
 $R(z_L;y)$& The loss function in the $L$-th layer.\\
 $\Omega_l(W_l)$& The regularization term in the $l$-th layer.\\
 $\varepsilon$& The tolerance of the nonlinear mapping.\\
\hline
 \end{tabular}
  \caption{ Notations used in this paper}
 \label{tab:notation}
 \end{table}
 \normalsize
\subsection{Inequality Approximation for Deep Learning}
\label{sec:problem}
\indent Important notations used in this paper are shown in Table \ref{tab:notation}. A typical MLP model consists of $L$ layers, each of which are defined by a linear mapping and a nonlinear activation function. A linear mapping is composed of  a weight vector $W_l\in \mathbb{R}^{n_l\times n_{l-1}}$ , where $n_l$ is the number of neurons on the $l$-th layer; a nonlinear mapping is defined by a  continuous activation function $h_l(\bullet)$. Given an input $a_{l-1}\in \mathbb{R}^{n_{l-1}}$ from the $(l-1)$-th layer, the $l$-th layer outputs $a_l=h_l(W_la_{l-1})$. By introducing an auxiliary variable $z_l$ as the output of the linear mapping, the neural network problem is formulated mathematically as follows:
 \begin{problem}
 \label{prob:original problem}
 \begin{align*}
       &\min\nolimits_{a_l,W_l,z_l} R(z_L;y)+\sum\nolimits_{l=1}^L\Omega_l(W_l),\\
     &s.t.\ z_l\!=\!W_l\!a_{l-1}\ (l\!=\!1,\!\cdots\!,\!L), \ a_l\!=\!h_l(z_l) \ (l\!=\!1,\!\cdots\!,\!L\!-\!1),
\end{align*}
 \end{problem}
where $a_0=x\in\mathbb{R}^d$ is the input of the neural network, $d$ is the number of feature dimensions, and $y$ is a predefined label vector. $R(z_L;y)\geq 0$ is a continuous loss function for the $L$-th layer, which is convex and proper,  and $\Omega_l(W_l)\geq 0$ is a regularization term on the $l$-th layer, which is also continuous, convex and proper.\\
\indent The equality constraint $a_l=h_l(z_l)$ is the most challenging one to handle here because common activation functions such as sigmoid \cite{glorot2010understanding} are nonlinear. This makes them nonconvex constraints and hence it is difficult to obtain the optimal solution when solving the $z_l$-subproblem \cite{taylor2016training}. Moreover, there is no guarantee for alternating minimization methods to solve the nonlinear equality constrained Problem \ref{prob:original problem} \cite{WANG2021100009}. To deal with these two challenges, the following assumption is required for problem transformation:
\begin{assumption}
\label{ass:quasilinear}
$h_l(z_l) \ (l=1,\dots,n)$ are quasilinear. 
\end{assumption}
The quasilinearity is defined in the appendix. Assumption \ref{ass:quasilinear} is so mild that most of the widely used nonlinear activation functions satisfy it, including tanh \cite{zamanlooy2014efficient}, sigmoid \cite{glorot2010understanding}, and the Rectified Linear Unit (ReLU) \cite{maas2013rectifier}. Then we innovatively transform the original nonconvex constraints into inequality constraints, which can be an infinite approximation of Problem  \ref{prob:original problem}. To do this, we  introduce a tolerance $\varepsilon>0$ and reformulate Problem  \ref{prob:original problem} to the following: 
\begin{align*}
     &\min\nolimits_{W_l,z_l,a_l} R(z_L;y)+\sum\nolimits_{l=1}^L\Omega_l(W_l),
      \\& s.t.\ z_l=W_la_{l-1}(l=1,\cdots,L),\\& h_l(z_l)-\varepsilon\leq a_l\leq h_l(z_l)+\varepsilon (l=1,\cdots,L-1).
\end{align*}
For the linear constraint $z_l=W_la_{l-1}$, this can be transformed into a penalty term in the objective function to minimize the difference between $z_l$ and $W_la_{l-1}$. The formulation is shown as follows:
\begin{problem}
\label{prob:inequality constrained DNN}
\begin{align*}
&\min\nolimits_{W_l,z_l,a_l} F(\textbf{W},\textbf{z},\textbf{a})\\
&=R(z_L;y)+\sum\nolimits_{l=1}^L\Omega_l(W_l)+\sum\nolimits_{l=1}^L\phi(a_{l-1},W_l,z_l),\\& s.t. \ h_l(z_l)-\varepsilon\leq a_l\leq h_l(z_l)+\varepsilon (l=1,\cdots,L-1).
\end{align*}
\end{problem}
   The penalty term is defined as $\phi(a_{l-1},W_l,z_l)=\frac{\rho}{2}\Vert z_l-W_la_{l-1}\Vert^2_2$, where $\rho>0$ a penalty parameter. $\textbf{W}=\{W_l\}_{l=1}^{L}$, $\textbf{z}=\{z_l\}_{l=1}^{L}$,$\textbf{a}=\{a_l\}_{l=1}^{L-1}$. As $\rho\rightarrow \infty$ and $\varepsilon\rightarrow 0$, Problem \ref{prob:inequality constrained DNN} approaches Problem  \ref{prob:original problem}.\\  
\indent The introduction of $\varepsilon$ is to project the nonconvex constraints to $\varepsilon$-balls, thus transforming the nonconvex Problem \ref{prob:original problem} into Problem \ref{prob:inequality constrained DNN}. Even though Problem \ref{prob:inequality constrained DNN} is still nonconvex because $h_l(z_l)$ can be nonconvex (e.g. tanh and smooth sigmoid), it is convex with regard to one variable when others are fixed (i.e. multi-convex), which is much easier to solve by alternating minimization \cite{xu2013block}. For example, Problem \ref{prob:inequality constrained DNN} is convex with regard to $\textbf{z}$ when $\textbf{W}$, and $\textbf{a}$ are fixed.  
\begin{algorithm} %算法开始
\caption{ The proposed mDLAM algorithm}
%算法的题目 
\small
\begin{algorithmic}[1]
\label{algo:mDLAM}
%此处的[1]控制一下算法中的每句前面都有标号 
\REQUIRE $y$, $a_0=x$. %输入条件(此处的REQUIRE默认关键字为Require) 
\ENSURE $a_l,W_l, z_l(l=1,\cdots,L)$. %输出结果(此处的ENSURE默认关键字为Ensure) 
\STATE Initialize $\rho$, $k=0$. $s^0=0$.
\REPEAT
\STATE $s^{k+1}\leftarrow \frac{1+\sqrt{1+4(s^k)^2}}{2}$
\STATE $\omega^k\leftarrow\frac{s^k-1}{s^{k+1}}$
\FOR{$l=1$ to $L$}
\STATE $\overline{W}^{k+1}_l\leftarrow {W}^{k}_l+({W}^{k}_l-{W}^{k-1}_l)\omega^k$ 
and update $W_l^{k+1}$ in Equation \eqref{eq:update W}.
\IF{$F(\textbf{W}^{k+1}_{\leq l},\textbf{z}^{k+1}_{\leq l-1},\textbf{a}^{k+1}_{\leq l-1})\geq F(\textbf{W}^{k+1}_{\leq l-1},\textbf{z}^{k+1}_{\leq l-1},\textbf{a}^{k+1}_{\leq l-1})$\hfill\COMMENT{\#$W_l^{k+1}$ increases the objective $F$}}
\STATE $\overline{W}^{k+1}_l\leftarrow {W}^{k}_l$ and update $W_l^{k+1}$ in Equation \eqref{eq:update W}.
\ENDIF
\STATE $\overline{z}^{k+1}_l\leftarrow {z}^{k}_l+({z}^{k}_l-{z}^{k-1}_l)\omega^k$
\IF{$l=L$}
\STATE Update $z_L^{k+1}$ in  Equation \eqref{eq:update zl}.
\IF{$F(\textbf{W}^{k+1}_{\leq L},\textbf{z}^{k+1}_{\leq L},\textbf{a}^{k+1}_{\leq L-1})\geq F(\textbf{W}^{k+1}_{\leq L},\textbf{z}^{k+1}_{\leq L-1},\textbf{a}^{k+1}_{\leq L-1})$\hfill\COMMENT{\#$z_L^{k+1}$ increases the objective $F$}}
\STATE $\overline{z}^{k+1}_L\leftarrow {z}^{k}_L$ and update $z_L^{k+1}$ in Equation \eqref{eq:update zl}.
\ENDIF
\ELSE
\STATE Update $z_l^{k+1}$ in  Equation \eqref{eq:update z}.
\IF{$F(\textbf{W}^{k+1}_{\leq l},\textbf{z}^{k+1}_{\leq l},\textbf{a}^{k+1}_{\leq l-1})\geq F(\textbf{W}^{k+1}_{\leq l},\textbf{z}^{k+1}_{\leq l-1},\textbf{a}^{k+1}_{\leq l-1})$\hfill\COMMENT{\#$z_l^{k+1}$ increases the objective $F$}}
\STATE $\overline{z}^{k+1}_l\leftarrow {z}^{k}_l$ and update $z_l^{k+1}$ in Equation \eqref{eq:update z}.
\ENDIF
\STATE $\overline{a}^{k+1}_l\leftarrow {a}^{k}_l+({a}^{k}_l-{a}^{k-1}_l)\omega^k$ and update $a^{k+1}_{l}$ in Equation \eqref{eq:update a}.
\IF{$F(\textbf{W}^{k+1}_{\leq l},\textbf{z}^{k+1}_{\leq l},\textbf{a}^{k+1}_{\leq l})\geq F(\textbf{W}^{k+1}_{\leq l},\textbf{z}^{k+1}_{\leq l},\textbf{a}^{k+1}_{\leq l-1})$\hfill\COMMENT{\#$a_l^{k+1}$ increases the objective $F$}}
\STATE $\overline{a}^{k+1}_l\leftarrow {a}^{k}_l$ and update $a_l^{k+1}$ in Equation \eqref{eq:update a}. 
\ENDIF
\ENDIF
\ENDFOR
\STATE $k\leftarrow k+1$.
\UNTIL convergence.
\STATE Output $a_l,W_l,z_l$.
\end{algorithmic}
\end{algorithm}
\subsection{Alternating Optimization}
\label{sec:quadratic approxmation}
\indent We present the mDLAM algorithm to solve Problem \ref{prob:inequality constrained DNN} in this section. A potential challenge to solve Problem \ref{prob:inequality constrained DNN} is a slow theoretical convergence rate. For example, the convergence rate of the dlADMM algorithm to solve Problem \ref{prob:inequality constrained DNN} is sublinear $o(1/k)$, where $k$ is the number of iterations \cite{wang2019admm}. In order to address this challenge, we apply the famous Nesterov acceleration technique to boost the convergence of our proposed mDLAM algorithm, and we prove its linear convergence theoretically in the next section.\\
\indent Algorithm \ref{algo:mDLAM} shows our proposed mDLAM algorithm. To simplify the notation, $\textbf{W}^{k+1}_{\leq l}=\{\{W^{k+1}_i\}_{i=1}^{l},\{W^{k}_i\}_{i=l+1}^{L}\}$, $\textbf{z}^{k+1}_{\leq l}=\{\{z^{k+1}_i\}_{i=1}^{l},\{z^{k}_i\}_{i=l+1}^{L}\}$ and $\textbf{a}^{k+1}_{\leq l}=\{\{a^{k+1}_i\}_{i=1}^{l},\{a^{k}_i\}_{i=l+1}^{L-1}\}$. In Algorithm \ref{algo:mDLAM}, Lines 6, 10, and 21 apply the Nestrov acceleration technique and update $W_l$, $z_l$ and $a_l$, respectively. the proposed mDLAM algorithm guarantees the decrease of objective $F$:  for example, if the updated $W^{k+1}_l$ in Line 7 of Algorithm \ref{algo:mDLAM} increases the value of $F$, i.e. $F(\textbf{W}^{k+1}_{\leq l},\textbf{z}^{k+1}_{\leq l-1},\textbf{a}^{k+1}_{\leq l-1})\geq F(\textbf{W}^{k+1}_{\leq l-1},\textbf{z}^{k+1}_{\leq l-1},\textbf{a}^{k+1}_{\leq l-1})$, then $W^{k+1}_l$ is updated again by setting $\overline{W}^{k+1}_l=W^{k}_l$ in Line 8 of Algorithm \ref{algo:mDLAM}, which ensures the decline of $F$. The same procedure is applied in Lines 13-15, Lines 18-20, and Lines 22-24 in Algorithm \ref{algo:mDLAM}, respectively.\\
\indent Next, all subproblems are shown as follows:\\
\textbf{1. Update $W_l$}\\
\indent The variables $W_l(l=1,\cdots,L)$ are updated as follows:
\begin{align}
    W^{k+1}_l\leftarrow\arg\min\nolimits_{W_l} \phi(a^{k+1}_{l-1},W_l,z^k_l)+\Omega_l(W_l).
    \label{eq: update W original}
\end{align}
\indent Because $W_l$ and $a_{l-1}$ are coupled in $\phi(\bullet)$, solving $W_l$ requires an inversion operation of $a^{k+1}_{l-1}$, which is computationally expensive. Motivated by the dlADMM algorithm \cite{wang2019admm},  we define $P^{k+1}_l(W_l;\theta^{k+1}_l)$ as a quadratic approximation of $\phi$ at $W^k_{l}$ as follows:
\begin{align*}&P^{k+1}_l(W_l;\theta^{k+1}_l)=\phi(a^{k+1}_{l-1},\overline{W}^{k+1}_{l},z^k_{l})\\&+(\nabla_{\overline{W}^{k+1}_{l}} \phi)^T(W_l-\overline{W}^{k+1}_{l})+\frac{\theta_l^{k+1}}{2}\Vert W_l-\overline{W}^{k+1}_{l}\Vert^2_2,
\end{align*}
where $\theta_l^{k+1}>0$ is a scalar parameter, which can be chosen by the backtracking algorithm \cite{wang2019admm} to meet the following condition 
\begin{align}
    P^{k+1}_l(W^{k+1}_l;\theta^{k+1}_l)\geq  \phi(a^{k+1}_{l-1},W^{k+1}_{l},z^k_{l}). \label{ineq:ineq1}
\end{align}Rather than minimizing  Equation \eqref{eq: update W original}, we instead minimize the following:
\begin{align}
 &W^{k+1}_l \leftarrow \arg\min\nolimits_{W_l} P^{k+1}_l(W_l;\theta_l^{k+1})+\Omega_l (W_l).     \label{eq:update W}
\end{align}
For  $\Omega_l(W_l)$, common  regularization terms like $\ell_1$ or $\ell_2$ regularizations lead to closed-form solutions.\\ 
\textbf{2. Update $z_l$}\\
\indent The variables $z_l(l=1,\cdots,L)$ are updated as follows:
\begin{align*}
      z^{k+1}_l &\leftarrow\arg\min\nolimits_{z_l} \phi(a^{k+1}_{l-1},W^{k+1}_l,z_l),\\& s.t. \ h_l(z_l)-\varepsilon\leq a_l\leq h_l(z_l)+\varepsilon \ (l< L).\\
    z^{k+1}_L &\leftarrow\arg\min\nolimits_{z_L} \phi(a^{k+1}_{L-1},W^{k+1}_L,z_L)+R(z_L;y).
\end{align*}
Similar to updating $W_l$, we define $V^{k+1}_l(z_l)$ as follows:
\begin{align*}
    &V^{k+1}_l(z_l)=\phi(a_{l-1}^{k+1},W^{k+1}_l,\overline{z}^{k+1}_{l})+(\nabla_{\overline{z}^{k+1}_{l}}\phi)^T(z_l-\overline{z}^{k+1}_{l})\\&+\frac{\rho}{2}\Vert z_l-\overline{z}^{k+1}_{l}\Vert^2_2.
\end{align*}
Hence, we solve the following problems:
\begin{align}
 \nonumber &z^{k+1}_l \leftarrow\arg\min\nolimits_{z_l} V^{k+1}_l(z_l), \\& s.t. \ h_l(z_l)-\varepsilon\leq a_l\leq h_l(z_l)+\varepsilon \ (l< L) \label{eq:update z}.\\
    z^{k+1}_L &\leftarrow\arg\min\nolimits_{z_L} V^{k+1}_L(z_L)+R(z_L;y). \label{eq:update zl}   
\end{align}
As for $z_l(l=1,\cdots,l-1)$, the solution is
\begin{align*}
z_l^{k+1}\leftarrow \min(\max(B^{k+1}_1,\overline{z}^{k+1}_{l}-\nabla\phi_{\overline{z}^{k+1}_{l}}/\rho), B^{k+1}_2),   
\end{align*} 
where $B^{k+1}_1$ and $B^{k+1}_2$ represent the lower bound and the upper bound of the set $\{ z_l|h_l(z_l)-\varepsilon\leq a^{k}_l\leq h_l(z_l)+\varepsilon\}$.   Equation \eqref{eq:update zl} is easy to solve using the Fast Iterative Soft Thresholding Algorithm (FISTA) \cite{beck2009fast}.\\
\textbf{3. Update $a_l$}\\
\indent The variables $a_l(l=1,\cdots,L-1)$ are updated as follows:
\begin{align*}
    a^{k+1}_{l}&\leftarrow \arg\min\nolimits_{a_{l}} \phi(a_{l},W^k_{l+1},z^k_{l+1}),\\& s.t. \ h_l(z^{k+1}_l)-\varepsilon\leq a_l\leq h_l(z^{k+1}_l)+\varepsilon.
\end{align*}
\indent Similar to updating $W^{k+1}_l$, $Q^{k+1}_l(a_{l};\tau^{k+1}_l)$ is defined as 
\begin{align*}
&Q^{k+1}_l(a_{l};\tau^{k+1}_l)=\phi(\overline{a}^{k+1}_{l},W^k_{l+1},z^k_{l+1})\\&+(\nabla_{\overline{a}^{k+1}_{l}} \phi)^T(a_{l}-\overline{a}^{k+1}_{l})+\frac{\tau_l^{k+1}}{2}\Vert a_{l}-\overline{a}^{k+1}_{l}\Vert^2_2,
\end{align*}
and this allows us to solve the following problem instead:
\begin{align} \nonumber a^{k+1}_l&\leftarrow\arg\min\nolimits_{a_l} Q^{k+1}_l(a_l;\tau^{k+1}_l),\\& s.t. \ h_l(z^{k+1}_l)-\varepsilon\leq a_l\leq h_l(z^{k+1}_l)+\varepsilon, \label{eq:update a}
\end{align}
where $\tau_l^{k+1}>0$ is a scalar parameter, which can be chosen by the backtracking algorithm \cite{wang2019admm} to meet the following condition:
\begin{align*}
 Q^{k+1}_l(a^{k+1}_{l};\tau^{k+1}_l)\geq \phi(a^{k+1}_{l},W^k_{l+1},z^k_{l+1}).   
\end{align*}
The solution can be obtained by
\begin{align*}
 a_{l}^{k+1}\!&\leftarrow\! \min(\max(h_l(z^{k+1}_l)\!-\!\varepsilon,\overline{a}^{k+1}_{l}\!-\!\nabla_{\overline{a}^{k+1}_{l}}\phi/\tau_{l}^{k+1}),\\&h_l(z^{k+1}_l)\!+\!\varepsilon).
 \end{align*}
 \section{Convergence Analysis}
\label{sec:convergence}
\indent In this section, the convergence of the proposed algorithm is analyzed. Due to space limit, all proofs are detailed in the appendix. The following mild assumption is required for the convergence analysis of the proposed mDLAM algorithm:
\begin{assumption}
$F(\textbf{W},\textbf{z},\textbf{a})$ is coercive over the domain $\{(\textbf{W},\textbf{z},\textbf{a})| h_l(z_l)-\varepsilon \leq a_l\leq h_l(z_l)+\varepsilon \ (l=1,\cdots,L-1)\}$.
\label{ass:coercive}
\end{assumption}
 The coercivity is defined in the Appendix. Assumption \ref{ass:coercive} is also mild such that common loss functions such as the least square loss and the cross-entropy loss satisfy it \cite{wang2019admm}. 
 \subsection{Convergence Properties}
 \indent Firstly, the following preliminary lemma is useful to prove the convergence properties of the proposed mDLAM algorithm.
 \begin{lemma}
\label{lemma:lemma 3}
In Algorithm \ref{algo:mDLAM}, there exist $\alpha_l^k, \gamma_l^k, \delta_l^k>0$ such that for $\forall k\in \mathbb{N}$, $W^k_l,z^k_l(l=1,2,\cdots,L)$, and $a^k_l(l=1,2,\cdots,L-1)$, it holds that
\begin{align}
\nonumber& F(\textbf{W}^{k+1}_{\leq l-1},\textbf{z}^{k+1}_{\leq l-1},\textbf{a}^{k+1}_{\leq l-1})-F(\textbf{W}^{k+1}_{\leq l},\textbf{z}^{k+1}_{\leq l-1},\textbf{a}^{k+1}_{\leq l-1})\\&\geq\frac{\alpha_l^{k+1}}{2}\Vert W^{k+1}_l-W^k_l\Vert^2_2,\label{eq:mDLAM w optimality}\\
\nonumber&F(\textbf{W}^{k+1}_{\leq l},\textbf{z}^{k+1}_{\leq l-1},\textbf{a}^{k+1}_{\leq l-1})-F(\textbf{W}^{k+1}_{\leq l},\textbf{z}^{k+1}_{\leq l},\textbf{a}^{k+1}_{\leq l-1})\\&\geq \frac{\gamma^{k+1}_l}{2}\Vert z^{k+1}_l-z^k_l\Vert^2_2, \label{eq:mDLAM z optimality}\\\nonumber&F(\textbf{W}^{k+1}_{\leq l},\textbf{z}^{k+1}_{\leq l},\textbf{a}^{k+1}_{\leq l-1})-F(\textbf{W}^{k+1}_{\leq l},\textbf{z}^{k+1}_{\leq l},\textbf{a}^{k+1}_{\leq l})\\& \geq \frac{\delta^{k+1}}{2}\Vert a^{k+1}_l-a^k_l\Vert^2_2.
\label{eq:mDLAM a optimality}
\end{align}
\end{lemma}
 \indent It shows that the objective decreases when all variables are updated. Based on Assumption \ref{ass:coercive} and Lemma \ref{lemma:lemma 3},  three convergence properties hold, which are shown in the following:
\begin{lemma}[Objective Decrease]
In Algorithm \ref{algo:mDLAM}, it holds that for any $k\in \mathbb{N}$,  $F(\textbf{W}^k,\textbf{z}^{k},\textbf{a}^{k})\geq F(\textbf{W}^{k+1},\textbf{z}^{k+1},\textbf{a}^{k+1})$. Moreover, $F$ is convergent. That is, $F(\textbf{W}^k,\textbf{z}^{k},\textbf{a}^{k})\rightarrow F^*$ as $k\rightarrow\infty$, where $F^*$ is the convergent value of $F$.
\label{lemma:objective decrease}
\end{lemma}
 This lemma guarantees the decrease and hence convergence of the objective. 
 \begin{lemma}[Bounded Objective and Variables]
In Algorithm \ref{algo:mDLAM}, it holds that for any $k\in \mathbb{N}$\\
(a). $\textbf{F}(\textbf{W}^k,\textbf{z}^k,\textbf{a}^k)$ is upper bounded. Moreover, $\lim\nolimits_{k\rightarrow\infty}\textbf{W}^{k+1}-\textbf{W}^{k}=0$, $\lim\nolimits_{k\rightarrow\infty}\textbf{z}^{k+1}-\textbf{z}^{k}=0$, and $\lim\nolimits_{k\rightarrow\infty}\textbf{a}^{k+1}-\textbf{a}^{k}=0$.\\
(b). $(\textbf{W}^k,\textbf{z}^k,\textbf{a}^k)$ is bounded. That is, there exist scalars $M_\textbf{W}$,$M_\textbf{z}$ and $M_\textbf{a}$ such that $\Vert\textbf{W}^k\Vert\leq M_\textbf{W}$, $\Vert\textbf{z}^k\Vert\leq M_\textbf{z}$ and $\Vert\textbf{a}^k\Vert\leq M_\textbf{a}$.
\label{lemma:bounded objective}
\end{lemma}
This lemma ensures that the objective and all variables are bounded in the proposed mDLAM algorithm. Moreover, the gap between the same variables in the neighboring iterations (e.g. $\textbf{W}^{k+1}$ and $\textbf{W}^{k}$) is convergent to $0$.
\begin{lemma}[Subgradient Bound]
\label{lemma:subgradient bound} 
In Algorithm \ref{algo:mDLAM},
there exist $C_2=\max(\rho M_{\textbf{a}},\rho M^2_{\textbf{a}}+\theta^{k+1}_1,\rho M^2_{\textbf{a}}+\theta^{k+1}_2,\cdots,\rho M^2_{\textbf{a}}+\theta^{k+1}_L)$, and $ g^{k+1}_1\in \partial_{\textbf{W} ^{k+1}}F$ such that for any $k\in \mathbb{N}$
\begin{align*}
    \Vert g^{k\!+\!1}_1\Vert\!\leq\!C_2(\Vert \textbf{W}^{k\!+\!1}\!-\!\textbf{W}^{k}\Vert\!+\!\Vert\textbf{z}^{k\!+\!1}\!-\!\textbf{z}^{k}\Vert\!+\!\Vert \textbf{W}^{k}\!-\!\textbf{W}^{k\!-\!1}\Vert).
\end{align*}
\end{lemma}
\indent The above lemma states that the subgradient of the objective is bounded by its variables. This suggests that the subgradient is convergent to $0$, and thus proves its convergence to a stationary point.
\subsection{Convergence of the proposed mDLAM Algorithm}
\label{sec:mDLAM convergence}
\indent Next we discuss the convergence of the proposed mDLAM algorithm. The first theorem  guarantees that the proposed mDLAM algorithm converges to a stationary point.
\begin{theorem}[Convergence to a Stationary Point]
\label{thero:mDLAM convergence} In Algorithm \ref{algo:mDLAM},
for $\textbf{W}$ in Problem \ref{prob:inequality constrained DNN}, for any $\rho>0$ and $\varepsilon>0$, starting from any $\textbf{W}^0$ , any limit point $\textbf{W}^*$ is a stationary point of Problem \ref{prob:inequality constrained DNN}. That is, $0\in \partial{_{\textbf{W}^*}} F$.
\end{theorem}
\indent As stated in Theorem \ref{thero:mDLAM convergence}, the convergence always holds no matter how $\textbf{W}$ is initialized, and whatever $\rho$ and $\varepsilon$ are chosen. It is better than the dlADMM algorithm \cite{wang2019admm}, which requires the hyperparameter to be sufficiently large.
\begin{theorem}[Linear Convergence Rate]
In Algorithm \ref{algo:mDLAM}, 
if $F$ is locally strongly convex, then for any $\rho$, there exist $\varepsilon>0$, $k_1\in \mathbb{N}$ and $0<C_1<1$ such that it holds for $k>k_1$  that
\begin{align*}
       F(\textbf{W}^{k\!+\!1},\textbf{z}^{k\!+\!1},\textbf{a}^{k+1})\!-\!F^*\!\leq\! C_1(F(\textbf{W}^{k\!-\!1},\textbf{z}^{k\!-\!1},\textbf{a}^{k\!-\!1})\!-\!F^*).
\end{align*}
\label{thero:mDLAM convergence rate}
\end{theorem}
Theorem \ref{thero:mDLAM convergence rate} shows that the proposed mDLAM algorithm converges linearly for sufficiently large iterations. Common loss functions like the square loss or the cross-entropy loss are locally strongly convex \cite{xu2013block}, which make $F$ locally strongly convex. Therefore, Theorem \ref{thero:mDLAM convergence rate} covers a wide range of loss functions. Compared with existing alternating minimization methods (e.g. dlADMM \cite{wang2019admm}) with a sublinear $o(1/k)$ convergence rate, the proposed mDLAM algorithm achieves a theoretically better linear convergence rate.
\subsection{Discussion}
\label{sec:discussion}
\indent We discuss convergence conditions of the proposed mDLAM algorithm compared with SGD-type methods and the dlADMM method. The comparison demonstrates that our convergence conditions are more general than others.\\
\textbf{1. mDLAM versus SGD}\\
\indent One influential work by Ghadimi et al. \cite{ghadimi2016accelerated} guaranteed that the SGD converges to a critical point, which is similar to our convergence results. While the SGD requires the objective function to be Lipschitz differentiable, bounded from below \cite{ghadimi2016accelerated}, our mDLAM allows for non-smooth functions such as ReLU. Therefore, our convergence conditions are milder than SGD.\\
\textbf{2. mDLAM versus dlADMM}\\
\indent Wang et al. \cite{wang2019admm} proposed an improved version of ADMM for deep learning models called dlADMM. They showed that the dlADMM is convergent to a critical point. However,  assumptions of our mDLAM are milder than those of the dlADMM: the mDLAM requires activation functions to be quasilinear, which includes sigmoid, tanh, ReLU, and leaky ReLU, while the dlADMM assumes that activation functions make subproblems solvable, which only includes ReLU and leaky ReLU. Such difference originates from different ways of addressing nonlinear activations: the dlADMM treats them as $L_2$ penalties. For tanh and sigmoid, subproblems are difficult to solve and may refer to lookup tables \cite{wang2019admm}. However, the mDLAM relaxes them via inequality constraints, and subproblems have closed-form solutions.\\
\section{Experiments}
\label{sec:experiment}
In this section, we evaluate the proposed mDLAM algorithm on four benchmark datasets. Convergence and efficiency are demonstrated. The performance of the proposed mDLAM algorithm is compared with several state-of-the-art optimizers. All experiments were conducted on a 64-bit machine with Intel(R) Xeon(R) Silver 4110 CPU and 64GB RAM. 
\begin{table}[]
    \centering
    \begin{tabular}{c|c|c|c|c}
    \hline\hline
         Dataset&\tabincell{c}{ Node\#}&\tabincell{c}{ Edge\#}&
\tabincell{c}{Class\#}&
\tabincell{c}{Feature\#}\\\hline Cora&2708& 5429&7&1433\\\hline Pubmed&19717&44338&3&500\\\hline Citeseer&3327&4732&6&3703\\\hline 
Coauthor CS &18333&81894&15&6805\\\hline\hline
    \end{tabular}
    \caption{Statistics of four benchmark datasets.}
    \label{tab:dataset}
\end{table}
\subsection{Datasets and Parameter Settings}
\begin{table}[]
    \centering
    \scriptsize
    \begin{tabular}{c|c|c|c|c|c}
    \hline\hline
         Method&\tabincell{c}{Hyper-\\ parameters}&Cora&Pubmed&Citeseer&\tabincell{c}{Coauthor\\ CS}  \\
         \hline
         mDLAM&$\rho$&$1\times 10^{-3}$&0.01&$5\times 10^{-3}$&$1\times 10^{-4}$\\\hline
        GD&$\alpha$&0.01&0.01&0.01&$5\times 10^{-3}$\\\hline
        Adadelta&$\alpha$&0.01&0.1&0.01&0.05\\\hline
        Adagrad&$\alpha$&$5\times 10^{-3}$&$5\times 10^{-3}$&0.01&$5\times 10^{-3}$\\\hline
        Adam&$\alpha$&$1\times 10^{-3}$&$5\times 10^{-4}$&$1\times 10^{-3}$&$1\times 10^{-3}$\\\hline
        dlADMM&$\rho$&$1\times 10^{-6}$&$1\times 10^{-6}$&$1\times 10^{-6}$&$1\times 10^{-6}$\\\hline\hline
    \end{tabular}
    \caption{Hyperparameter settings on four datasets: they were chosen based on training performance.}
    \label{tab:hyperparameter}
\end{table}
\indent An important application of the MLP model is node classification on a graph based on augmented node features \cite{chen2021graph}. Specifically, given an adjacency matrix $A$ and a node feature matrix $H$ of a graph, we let the $k$-th augmented feature $X^k=HA^k(k=0,1,\cdots,4)$, which encodes information of graph topology via $A^k$, and then concatenate them into the input $X=[X_0,\cdots,X_4]$ \cite{chen2021graph}.  The MLP model is used to predict the node class based on the input $X$. We set up an architecture of three layers, each of which has 100 hidden units. The activation function was set to ReLU. The number of epoch was set to $200$. We test our model on four benchmark datasets: Cora \cite{sen2008collective}, Pubmed \cite{sen2008collective}, Citeseer \cite{sen2008collective} and Coauthor CS \cite{shchur2018pitfalls}, whose statistics are shown in Table \ref{tab:dataset}.\\
\indent Gradient Descent (GD) \cite{bottou2010large}, Adaptive learning rate method (Adadelta) \cite{zeiler2012adadelta}, Adaptive gradient algorithm (Adagrad) \cite{duchi2011adaptive}, Adaptive momentum estimation (Adam) \cite{kingma2014adam}, and deep learning Alternating Direction Method of Multipliers (dlADMM) \cite{wang2019admm} are state-of-the-art methods and hence were served as comparison methods. The full batch dataset was used for training models. All parameters were chosen by maximizing the accuracy of training datasets. Table \ref{tab:hyperparameter} shows hyperparameters of all methods: for the proposed mDLAM algorithm, $\rho$ controls quadratic terms in Problem \ref{prob:inequality constrained DNN}; $\alpha$ is a learning rate in the comparison methods except for dlADMM. $\rho$ controls a linear constraint in the dlADMM algorithm. The other hyperparameter $\varepsilon$ is chosen adaptively as follows: $\varepsilon^{k+1}=\max(\varepsilon^k/2,0.001)$ with $\varepsilon^0=100$. This makes inequality constraints relaxed at the early stage (i.e. $\varepsilon^k$ is large and hence constraints are easy to satisfy) and then tightens them as the mDLAM iterates.
\begin{figure}
    \centering
    \includegraphics[width=\linewidth]{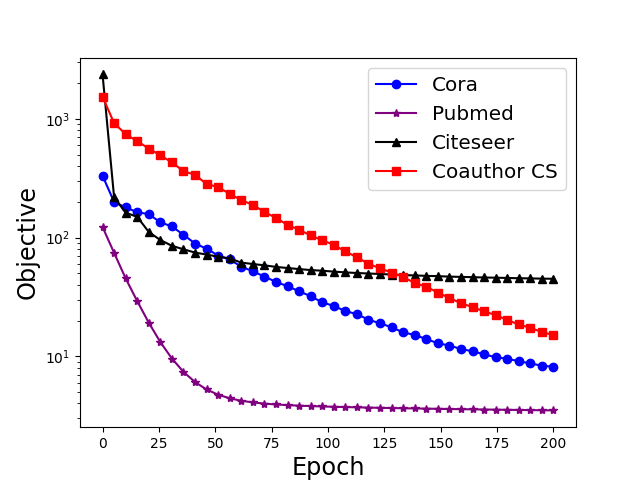}
    \caption{Convergence curves on four datasets: they all converge linearly when the epoch is larger than 100.}
    \label{fig:convergence}
\end{figure}
\subsection{Convergence}
Firstly, we investigate the convergence of the proposed mDLAM algorithm on four benchmark datasets using the hyperparameters summarized in Table \ref{tab:hyperparameter}. The relationship between the objective and the number of epochs is shown in Figure \ref{fig:convergence}. Overall, the objectives on the four datasets all decrease monotonically, which demonstrates the convergence of the proposed mDLAM algorithm. Nevertheless, objective curves vary in tendency: the curves on the Cora and Pubmed datasets drop drastically at the beginning and then reach the plateau when the epoch is around 75, while the curves on the other two datasets keep a downward tendency in the entire 200 epochs. Moreover, the objective on the Pubmed dataset is the lowest at the end of the training, while the objective on the Citeseer dataset is in the vicinity of 80, at least $60\%$ higher than objectives on the remaining datasets. It is easy to observe that all curves decline linearly when the epoch is higher than 100. This validates the linear convergence rate of our proposed mDLAM algorithm (i.e. Theorem \ref{thero:mDLAM convergence rate}).
\begin{figure}
    \centering
     \begin{minipage}{0.49\linewidth}
    \includegraphics[width=\linewidth]{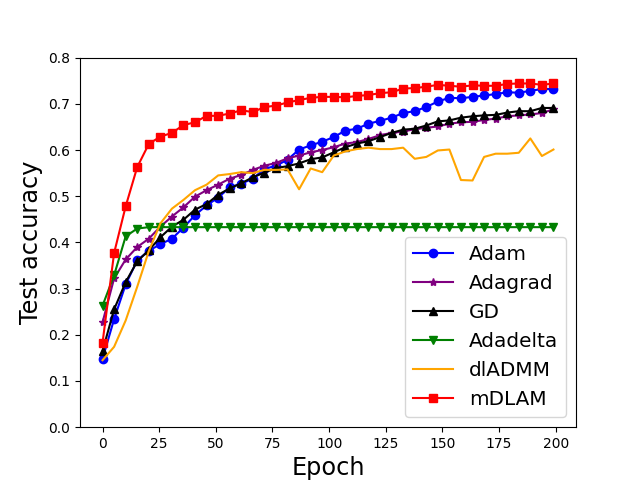}
    \centerline{(a). Cora.}
    \end{minipage}
    \hfill
        \begin{minipage}{0.49\linewidth}
    \includegraphics[width=\linewidth]{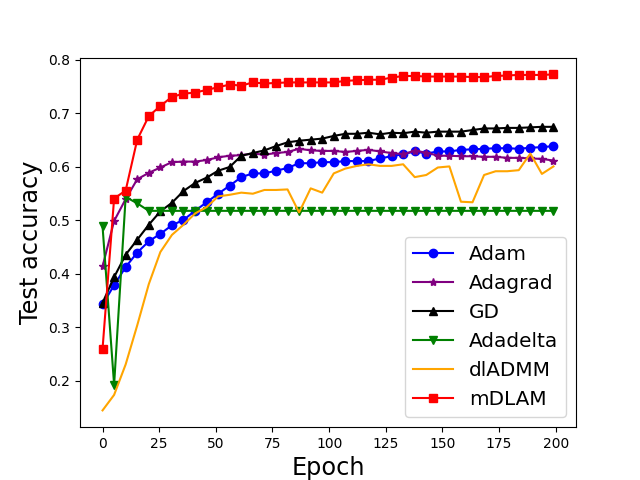}
    \centerline{(b). Pubmed.}
    \end{minipage}
           \vfill 
           \begin{minipage}{0.49\linewidth}
    \includegraphics[width=\linewidth]{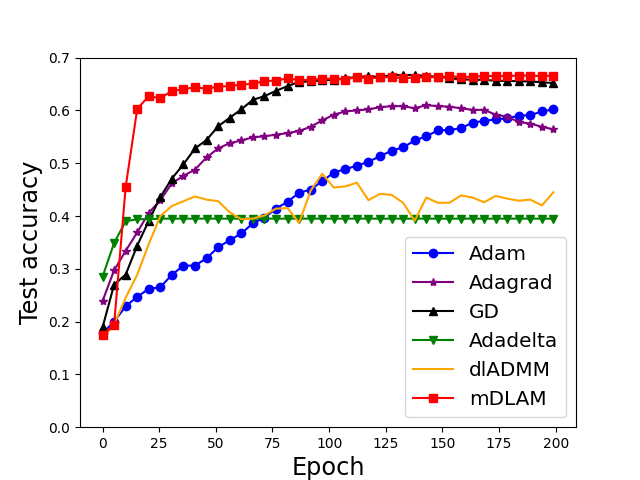}
    \centerline{(c). Citeseer.}
    \end{minipage}
               \begin{minipage}{0.49\linewidth}
    \includegraphics[width=\linewidth]{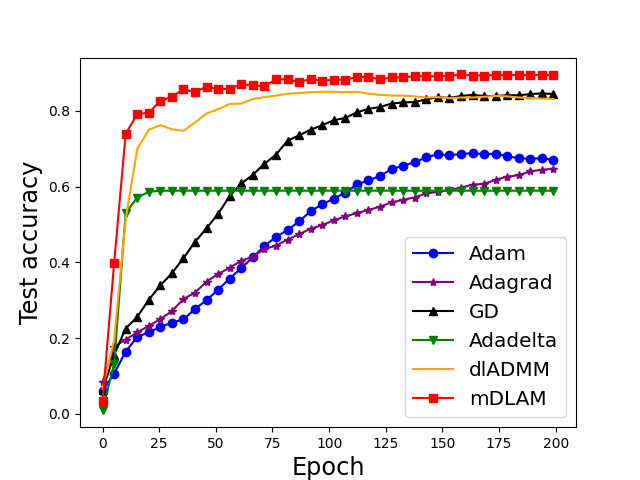}
    \centerline{(d). Coauthor CS.}
    \end{minipage}
    \caption{Test accuracy of all methods: the proposed mDLAM algorithm outperforms all other comparison methods in four datasets.}
    \label{fig:test performance}
\end{figure}

\subsection{Performance}
Next, the performance of the proposed mDLAM algorithm is compared against five state-of-the-art methods, as is illustrated in Figure \ref{fig:test performance}. X-axis and Y-axis represent epoch and test accuracy, respectively. Overall, the proposed mDLAM algorithm is superior to all other algorithms on four datasets, which has not only the highest test accuracy but also the fastest convergence speed. For example, the proposed mDLAM achieves $70\%$ test accuracy on the Cora dataset when the epoch is 100, while GD only attains $60\%$, and the Adadelta reaches the plateau of around $40\%$; As another example, the test accuracy of the proposed mDLAM on the Coauthor CS dataset is over $80\%$ at the $25$-th epoch, whereas most comparison methods such as Adam and GD reach half of its accuracy (i.e. $40\%$).  The Adadelta algorithm performs the worst among all comparison methods: it converges to a low test accuracy at the early stage, which is usually half of the accuracy accomplished by the proposed mDLAM algorithm. The other four comparison methods except Adagrad are on par with mDLAM in some cases: for example, the curves of dlADMM and GD are marginally behind that of mDLAM on the Coauthor CS dataset, and the performance of Adam almost reaches that of mDLAM on the Cora dataset. It is interesting to observe that curves of some methods decline at the end of 200 epochs such as the Adagrad on the Pubmed dataset and the Adam on the Coauthor CS dataset.
\begin{figure*}
    \centering
    \begin{minipage}{0.48\linewidth}
            \includegraphics[width=\linewidth]{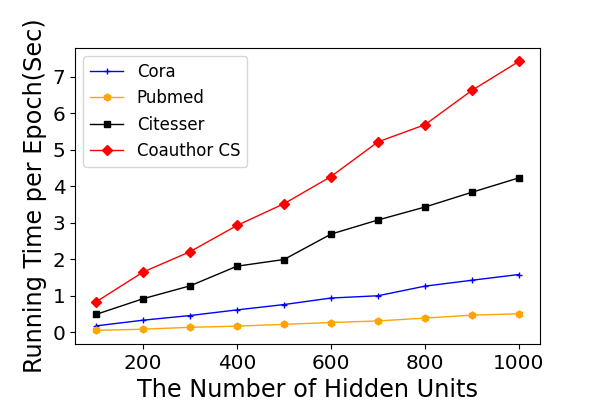}
            \centerline{(a). Running time versus the number of hidden units.}
    \end{minipage}
    \hfill
     \begin{minipage}{0.48\linewidth}
            \includegraphics[width=\linewidth]{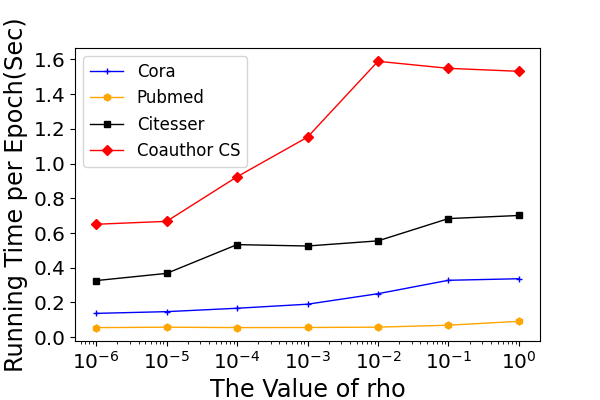}
            \centerline{(b). Running time versus the value of $\rho$.}
    \end{minipage}
    \caption{The relationship between the running time and: (a) the number of hidden units; (b) the value of $\rho$: the running time increases linearly with them in general.}
    \label{fig:running time}
\end{figure*}
\subsection{Sensitivity Analysis}
\indent We explore concerning factors of the running time and the test accuracy in this section.
\subsubsection{Running Time}
Moreover, it is important to explore the running time of the proposed mDLAM concerning two factors: the number of hidden units and the value of $\rho$. The running time was averaged by 200 epochs. Figure \ref{fig:running time}(a) depicts the relationship between the running time and the number of hidden units on four datasets, where the number of hidden units ranges from $100$ to $1,000$. The running times on all datasets are below 1 second per epoch when the number of hidden units is 100, and increase linearly with the number of hidden units in general. However, the rates of increase vary on different datasets: the curve on the Coauthor CS dataset has the sharpest slope, which reaches seven seconds per epoch when 1000 hidden units are applied, while the curve on the Pubmed dataset climbs slowly, which never surpasses 1 second. The curves on the Cora and the Citeseer datasets demonstrate a steady increase. \\
\indent To investigate the relationship between the running time per epoch and the value of $\rho$, we change $\rho$ from $10^{-6}$ to 1 while fixing others. Similar to Figure \ref{fig:running time}(a), the running time per epoch demonstrates a linear increase concerning the value of $\rho$ in general, as shown in Figure \ref{fig:running time}(b). Specifically, the curve on the Coauthor CS dataset is still the highest in slope, whereas the slope on the Pubmed dataset is the lowest. Moreover, the effect of the value of $\rho$ is less obvious than the number of hidden units. For example, in Figure \ref{fig:running time}(b) when $\rho$ is enlarged from $10^{-6}$ to $10^{-2}$, the running time on the Coauthor CS dataset merely ascends from around 0.65 to  1.6, while the increment of the running time on other datasets is less than 0.2. Moreover, a larger $\rho$ may reduce the running time. For instance, when $\rho$ increases from $10^{-2}$ to $1$, the running time on the Coauthor CS dataset drops slightly from 1.6 seconds to 1.5 seconds per epoch. The running times on the Cora and the Citeseer datasets climb steadily.
\begin{table}[]
    \centering
    \begin{tabular}{c|c|c|c|c|c}
    \hline\hline
    \multicolumn{6}{c}{Cora}\\\hline
         Epoch&40&80&120&160&200  \\
         \hline
         $\rho=1\times 10^{-4}$&0.677&0.695&0.695&0.693&0.692\\\hline
         $\rho=1\times 10^{-3}$&0.664&0.701&0.721&0.737&0.742\\\hline
         $\rho=1\times 10^{-2}$&0.562&0.581&0.604&0.623&0.638\\\hline
         \multicolumn{6}{c}{Pubmed}\\\hline
         Epoch&40&80&120&160&200  \\\hline
         $\rho=1\times 10^{-4}$&0.471&0.407&0.407&0.407&0.407\\\hline
         $\rho=1\times 10^{-3}$&0.663&0.645&0.640&0.650&0.649\\\hline
         $\rho=1\times 10^{-2}$&0.743&0.758&0.762&0.768&0.773\\\hline
         \multicolumn{6}{c}{Citeseer}\\\hline
         Epoch&40&80&120&160&200  \\
         \hline
         $\rho=1\times 10^{-4}$&0.528&0.529&0.530&0.531&0.535\\\hline
         $\rho=1\times 10^{-3}$&0.651&0.665&0.664&0.664&0.666\\\hline
         $\rho=1\times 10^{-2}$&0.631&0.638&0.642&0.648&0.653\\\hline
          \multicolumn{6}{c}{Coauthor CS}\\\hline
         Epoch&40&80&120&160&200  \\
         \hline
         $\rho=1\times 10^{-4}$&0.843&0.881&0.888&0.896&0.894\\\hline
         $\rho=1\times 10^{-3}$&0.780&0.807&0.825&0.839&0.835\\\hline
         $\rho=1\times 10^{-2}$&0.688&0.719&0.724&0.737&0.738\\
         \hline\hline
    \end{tabular}
    \caption{The effect of $\rho$ on test accuracy on four datasets: it affects performance significantly.}
    \label{tab:rho performance}
\end{table}
\begin{table}[]
    \centering
    \begin{tabular}{c|c|c|c|c|c}
    \hline\hline
    \multicolumn{6}{c}{Cora}\\\hline
         Epoch&40&80&120&160&200  \\
         \hline
         $\varepsilon^0=1$&0.620&0.679&0.712&0.735&0.743\\\hline
         $\varepsilon^0=10$&0.646&0.689&0.718&0.741&0.741\\\hline
         $\varepsilon^0=100$&0.664&0.701&0.721&0.737&0.742\\\hline
         \multicolumn{6}{c}{Pubmed}\\\hline
         Epoch&40&80&120&160&200  \\
         \hline
         $\varepsilon^0=1$&0.717&0.744&0.756&0.759&0.763\\\hline
         $\varepsilon^0=10$&0.731&0.753&0.759&0.762&0.765\\\hline
         $\varepsilon^0=100$&0.743&0.758&0.762&0.768&0.773\\\hline
         \multicolumn{6}{c}{Citeseer}\\\hline
         Epoch&40&80&120&160&200  \\
         \hline
         $\varepsilon^0=1$&0.564&0.615&0.638&0.653&0.663\\\hline
         $\varepsilon^0=10$&0.584&0.626&0.643&0.657&0.662\\\hline
         $\varepsilon^0=100$&0.640&0.656&0.664&0.663&0.668\\\hline
          \multicolumn{6}{c}{Coauthor CS}\\\hline
         Epoch&40&80&120&160&200  \\
         \hline
         $\varepsilon^0=1$&0.834&0.875&0.887&0.893&0.894\\\hline
         $\varepsilon^0=10$&0.852&0.866&0.892&0.893&0.893\\\hline
         $\varepsilon^0=100$&0.843&0.881&0.888&0.896&0.894\\
         \hline\hline
    \end{tabular}
    \caption{The effect of the initial value of $\varepsilon$ (i.e. $\varepsilon^0$) on test accuracy on four datasets: it only affects the convergence speed, but have little effect on final performance.}
    \label{tab:espilon performance}
\end{table}
\subsubsection{Test Accuracy}
\indent Finally, we investigate the effects of hyperparameters on test accuracy, namely, the value of $\rho$ and $\varepsilon$. Because $\varepsilon$ is dynamically set, we test its initial value $\varepsilon^0$. Table \ref{tab:rho performance} demonstrates the relationship between test accuracy and $\rho$ on four datasets. $\rho$ was chosen from $\{1\times 10^{-4},1\times 10^{-3},1\times 10^{-2}\}$. Overall, the choice of $\rho$ has a significant effect on the test accuracy. For example, when $\rho$ is changed from $1\times 10^{-4}$ to $1\times 10^{-3}$ on the Pubmed dataset, the performance has improved by approximately $60\%$, and the gain of performance is even roughly $90\%$ if it is modified to $1\times 10^{-2}$. On other datasets, the change of $\rho$ affects test accuracy by around $20\%$. For instance, the test accuracy on the Cora dataset and the Coauthor CS dataset can be improved to 0.74 and 0.89 if we set $\rho=1\times 10^{-3}$ and $\rho=1\times 10^{-4}$, respectively. The test accuracy on the Citeseer dataset is relatively robust to the change of $\rho$. As $\rho$ varies from $1\times 10^{-3}$ to $1\times 10^{-2}$, the test accuracy remains stable. Obviously, the test accuracy generally increases as the proposed mDLAM algorithm iterates. However, there are some exceptions: for example, the test accuracy has dropped slightly from $0.66$ to $0.65$ when $\rho=1\times 10^{-3}$ on the Pubmed dataset.\\
\indent Table \ref{tab:espilon performance} shows the relationship between test accuracy and the initial value of $\varepsilon$ (i.e. $\varepsilon^0$) on four datasets. $\varepsilon^0$ was chosen from $\{1,10,100\}$. It is obvious that test accuracy is resistant to the change of $\varepsilon^0$. For example, the test accuracy on the Coauthor CS dataset is in the vicinity of $0.89$ no matter whatever $\varepsilon$ is chosen. Moreover, the larger a $\varepsilon^0$ is, the faster convergence speed the proposed mDLAM algorithm gains. For instance, when $\varepsilon=100$, the test accuracy is 0.08 better than that in the case where $\varepsilon=1$ on the Citeseer dataset. Compared with Tables \ref{tab:rho performance} and \ref{tab:espilon performance}, the effect of $\rho$ is more signifcant than that of $\varepsilon^0$.
\section{Conclusion}
\label{sec:conclusion}
\indent  In this paper, we propose a novel formulation of the original neural network problem and a novel monotonous Deep Learning Alternating Minimization (mDLAM) algorithm. Specifically, the nonlinear constraint is projected into a convex set so that all subproblems are solvable. The Nesterov acceleration technique is applied to boost the convergence of the proposed mDLAM algorithm. Furthermore, a mild assumption is established to prove the convergence of our mDLAM algorithm. Our mDLAM algorithm can achieve a linear convergence rate, which is theoretically better than existing alternating minimization methods. The effectiveness of the proposed mDLAM algorithm is demonstrated via the outstanding performance on four benchmark datasets compared with state-of-the-art optimizers.
\section*{Acknowledgement}
This work was supported by the National Science Foundation (NSF) Grant No. 1755850, No. 1841520, No. 2007716, No. 2007976, No. 1942594, No. 1907805, a Jeffress Memorial Trust Award, Amazon Research Award, NVIDIA GPU Grant, and Design Knowledge Company (subcontract No: 10827.002.120.04).
\bibliographystyle{plain}
\bibliography{example_paper}
\normalsize
\newpage
\onecolumn
\appendix
\section*{Appendix}
\label{sec: aop}
\section{Definition}
\label{sec: definition}
Several definitions are shown here for the sake of convergence analysis.
\begin{definition}[Coercivity]
Any arbitrary function $G_2(x)$ is coercive over a nonempty set $dom(G_2)$ if as $\Vert x\Vert\rightarrow\infty$ and $x\in dom(G_2)$, we have $G_2(x)\rightarrow \infty$, where $dom(G_2)$ is a domain set of $G_2$.
\end{definition}
{\begin{definition}[Multi-convexity]
A function $f(x_1,x_2,\cdots,x_m)$ is a multi-convex function if $f$ is convex with regard to $x_i(i=1,\cdots,m)$ while fixing other variables.
\end{definition}
\begin{definition}[Lipschitz Differentiability]
A function $f(x)$ is Lipschitz differentiable with Lipschitz coefficient $L>0$ if for any $x_1,x_2\in \mathbb{R}$, the following inequality holds:
\begin{align*}
    \Vert\nabla f(x_1)-\nabla f(x_2)\Vert \leq L\Vert x_1-x_2\Vert.
\end{align*}
\end{definition}
For Lipschitz differentiability, we have the following lemma (Lemma 2.1 in \cite{beck2009fast}):
\begin{lemma}
\label{lemma:lipschitz differentiable}
If $f(x)$ is Lipschitz differentiable with $L>0$, then for any $x_1,x_2\in \mathbb{R}$
\begin{align*}
   f(x_1)\leq f(x_2)+ \nabla f^T(x_2)(x_1-x_2)+\frac{L}{2}\Vert x_1-x_2\Vert^2.
\end{align*}
\end{lemma}
\begin{definition}[Fr\'echet Subdifferential]
For each $x_1\in dom(u_1)$, the Fr\'echet subdifferential of $u_1$ at $x_1$, which is denoted as $\hat\partial u_1(x_1)$, is the set of vectors $v$, which satisfy
\begin{align*}
\lim_{x_2\neq x_1}\inf_{x_2\rightarrow x_1} (u_1(x_2)-u_1(x_1)-v^T(x_2-x_1))/\Vert x_2-x_1\Vert\geq 0.   
\end{align*}
 The vector $v\in \hat\partial u_1(x_1)$ is a Fr\'echet subgradient.
\end{definition}
\indent Then the definition of the limiting subdifferential, which is based on Fr\'echet subdifferential, is given in the following \cite{rockafellar2009variational}:
\begin{definition}[Limiting Subdifferential]
\label{def: limiting subdifferential}
For each $x\in dom(u_2)$, the limiting subdifferential (or subdifferential) of $u_2$ at $x$ is 
\begin{align*}
    \partial u_2(x)=&\{v_1| \exists \ x^k \rightarrow x, s.t. \ u_2(x^k) \rightarrow u_2(x), v^k\in\hat\partial u_2({x^k}), v^k\rightarrow v\}.
\end{align*}
where $x^k$ is a sequence whose limit is $x$ and the limit of $u_2(x^k)$ is $u_2(x)$, $v^k$ is a sequence, which is a Fr\'echet subgradient of $u_2$ at $x^k$ and whose limit is $v$.
The vector $v\in \partial u_2(x)$ is a limiting subgradient.
\end{definition}
\indent Specifically, when $u_2$ is convex, its limiting subdifferential is reduced to regular subdifferential \cite{rockafellar2009variational}, which is defined as follows:
\begin{definition}[Regular Subdifferential] For each $x_1\in dom(f)$, the regular subdifferential of a convex function $f$ at $x_1$, which is denoted as $\partial f(x_1)$, is the set of vectors $v$, which satisfy
\begin{align*}
f(x_2)\geq f(x_1)+v^T(x_2-x_1).
\end{align*}
 The vector $v\in \partial f(x_1)$ is a regular subgradient.
\end{definition}
\begin{definition}[Quasilinearity]
A function $f(x)$ is quasiconvex if for any sublevel set $S_\nu(f)=\{x|f(x)\leq \nu\}$ is a convex set. Likewise,  A function $f(x)$ is quasiconcave if for any superlevel set $S_\nu(f)=\{x|f(x)\geq \nu\}$ is a convex set. A function $f(x)$ is quasilinear if it is both quasiconvex and quasiconcave.
\end{definition}
\begin{definition}[Locally Strong Convexity]
A function $f(x)$ is locally strongly convex within a bound set $\mathbb{D}$ with a constant $\mu$ if
\begin{align*}
    f(y)\geq f(x)+g^T(y-x)+\frac{\mu}{2}\Vert x-y\Vert^2_2 \ \forall \ g\in \partial f(x) \text{ and } \ x,y\in \mathbb{D}.
\end{align*}
\end{definition}
Simply speaking, a locally strongly convex function lies above a quadratic function within a bounded set. 

}

\begin{definition}[Kurdyka-Lojasiewicz (KL) Property]
A function $f(x)$ has the KL Property at $\overline{x} \in dom \ \partial f= \{x\in \mathbb{R}: \partial f(x)\neq \emptyset \}$ if there exists $\eta\in (0,+\infty]$, a neighborhood $X$ of $\overline{x}$ and a function $\psi\in \Psi_\eta$, such that for all
\begin{align*}
    x\in X\cap \{x\in \mathbb{R}: f(\overline{x})<f(x)<f(\overline{x})+\eta\},
\end{align*}
the following inequality holds
\begin{align*}
    \psi^{'}(f(x)-f(\overline{x}))dist(0,\partial f(x))\geq 1,
    \end{align*}
    where $\Psi_\eta$ stands for a class of function $\psi: [0,\eta]\rightarrow \mathbb{R}^{+}$ satisfying: (1). $\phi$ is concave and $\psi^{'}(x)$ continuous on $(0,\eta)$; (2). $\psi$ is continuous at 0, $\psi(0)=0$; and (3). $\psi^{'}(x)>0, \forall x\in(0,\eta)$.
\end{definition}
The following lemma shows that a locally strongly convex function satisfies the KL Property:
\begin{lemma}[\cite{xu2013block}]
\label{lemma:locally strongly convex KL property}
A locally strongly convex function $f(x)$ with a constant $\mu$ satisfies the KL Property  at any $x\in \mathbb{D}$ with $\psi(x)=\frac{2}{\mu}\sqrt{x}$ and $X=\mathbb{D}\cap \{y:f(y)\geq f(x)\}$.
\end{lemma}
\section{Preliminary Results}
\label{sec:preliminary}
\indent In this section, we present preliminary lemmas of the proposed mDLAM algorithm. The limiting subdifferential is used to prove the convergence of the proposed mDLAM algorithm in the following convergence analysis. Without loss of generality,  $\partial R$ and $\partial \Omega_l(l=1,\cdots,n)$  are assumed to be  nonempty, and the limiting subdifferential of $F$ defined in Problem \ref{prob:inequality constrained DNN} is \cite{xu2013block}:\\
\begin{align*}
    \partial  F(\textbf{W},\textbf{z},\textbf{a})=\partial_{\text{W}} F\times \partial_{\text{z}} F \times \partial_{\text{a}} F,
\end{align*}
where $\times$ means the Cartesian product.
\begin{lemma}
\label{lemma:lemma 1}
If Equation \eqref{eq:update W} holds, then there exists $p\in \partial\Omega_l(W^{k+1}_l)$, the subgradient of $\Omega_l(W^{k+1}_l)$  such that
\begin{align*}
    \nabla_{\overline{W}^{k+1}_l} \phi +\theta^{k+1}_l (W^{k+1}_l-\overline{W}^{k+1}_l)+p=0.
\end{align*}
Likewise,  if Equation \eqref{eq:update z} holds, then there exists $q$ such that
\begin{align*}
    \nabla_{\overline{z}^{k+1}_l} \phi +\rho(z^{k+1}_l-\overline{z}^{k+1}_l)+q=0,
\end{align*}
where $q$ is a subgradient with regard to $z^{k+1}_l$ to satisfy the constraint $h_l(z^{k+1}_l)-\varepsilon\leq a^{k}_l\leq h_l(z^{k+1}_l)+\varepsilon$.
 If Equation \eqref{eq:update zl} holds, then there exists $u\in\partial R(z^{k+1}_L;y)$ such that
\begin{align*}
  \nabla_{\overline{z}^{k+1}_L} \phi +\rho(z^{k+1}_L-\overline{z}^{k+1}_L)+u=0.  
\end{align*}
 If Equation \eqref{eq:update a} holds, then there exists $v$ such that
\begin{align*}
 \nabla_{\overline{a}^{k+1}_l} \phi +\tau^{k+1}_l(a^{k+1}_l-\overline{a}^{k+1}_l)+v=0,  
\end{align*}
where $v$ is a subgradient with regard to $a^{k+1}_l$ to satisfy the constraint $h_l(z^{k+1}_l)-\varepsilon\leq a^{k+1}_l\leq h_l(z^{k+1}_l)+\varepsilon$.
\end{lemma}
\begin{proof}
    These can be obtained by directly applying the optimality conditions of  Equation \eqref{eq:update W},  Equation \eqref{eq:update z},  Equation \eqref{eq:update zl} and  Equation \eqref{eq:update a}, respectively.
\end{proof}
\begin{lemma}
For  Equation \eqref{eq:update z} and  Equation \eqref{eq:update zl}, the following inequalities hold:
\begin{align}
    &V^{k+1}_{l}(z^{k+1}_l)\geq \phi(a^{k+1}_{l-1},W^{k+1}_l,z^{k+1}_l).
    \label{eq:lipschitz z}
\end{align} 
\label{lemma:lemma 2}
\end{lemma}
\begin{proof}
Because $\phi(a_{l-1},W_l,z_l)$ is Lipschitz differentiable  with respect to $z_l$ with Lipschitz
coefficient $\rho$, we directly apply Lemma \ref{lemma:lipschitz differentiable} to $\phi$ to obtain Equation \eqref{eq:lipschitz z}.
\end{proof}
 \section{Main Proofs}
  \label{sec:main proofs}
\textbf{Proof of Lemma \ref{lemma:lemma 3}}
\begin{proof}
In Algorithm \ref{algo:mDLAM}, we only show Equation \eqref{eq:mDLAM w optimality}  because    Equation \eqref{eq:mDLAM z optimality} and  Equation \eqref{eq:mDLAM a optimality} follow the same routine of  Equation \eqref{eq:mDLAM w optimality}.\\
\indent In Line 7 of Algorithm \ref{algo:mDLAM}, if $ F(\textbf{W}^{k+1}_{\leq l},\textbf{z}^{k+1}_{\leq l-1},\textbf{a}^{k+1}_{\leq l-1})<F(\textbf{W}^{k+1}_{\leq l-1},\textbf{z}^{k+1}_{\leq l-1},\textbf{a}^{k+1}_{\leq l-1})$, then obviously there exists $\alpha^{k+1}_l>0$ such that Equation \eqref{eq:mDLAM w optimality} holds. Otherwise, according to Line 8 of Algorithm \ref{algo:mDLAM}, because $\Omega_{W_l}(W_l)$ and $\phi(a_{l-1},W_l,z_l)$ are convex with regard to $W_l$, according to the definition of regular subgradient, we have
\begin{align}
  &\Omega_l(W^k_{l})\geq  \Omega_l(W^{k+1}_l)+p^T(W_l^{k}-W_l^{k+1}) \label{ineq:ineq2}\\
  &\phi(a^{k+1}_{l-1},W^{k}_l,z^k_l)\geq  \phi(a^{k+1}_{l-1},\overline{W}^{k+1}_l,z^k_l)+\nabla_{\overline{W}^{k+1}_l}\phi^T(W_l^{k}-\overline{W}_l^{k+1}),
  \label{ineq:ineq4}
\end{align}
where $p$ is defined in the premise of Lemma \ref{lemma:lemma 1}. Therefore, we have 
\begin{align*}
    \nonumber&F(\textbf{W}^{k+1}_{\leq l-1},\textbf{z}^{k+1}_{\leq l-1},\textbf{a}^{k+1}_{\leq l-1})-F(\textbf{W}^{k+1}_{\leq l},\textbf{z}^{k+1}_{\leq l-1},\textbf{a}^{k+1}_{\leq l-1})\\\nonumber&=\phi(a^{k+1}_{l-1},W^{k}_l,z^k_l)+\Omega_l(W^k_l)-\phi(a^{k+1}_{l-1},W^{k+1}_l,z^k_l)-\Omega_l(W^{k+1}_{l})\  \text{(Definition of $F$ in Problem \ref{prob:inequality constrained DNN})}\\\nonumber&\geq \Omega_l(W^k_l)-\Omega_l(W^{k+1}_{l})-(\nabla_{\overline{W}^{k+1}_l}\phi)^T(W^{k+1}_l-\overline{W}^{k+1}_l)-\frac{\theta_l^{k+1}}{2}\Vert W^{k+1}_l-\overline{W}^{k+1}_l\Vert^2_2-\phi(a^{k+1}_{l-1},\overline{W}^{k+1}_l,z^k_l)\\\nonumber&+\phi(a^{k+1}_{l-1},W^{k}_l,z^k_l)(\text{Equation \eqref{ineq:ineq1}})\\\nonumber&\geq p^T(W^{k}_l-W^{k+1}_l)-(\nabla_{\overline{W}^{k+1}_l}\phi)^T(W^{k+1}_l-{W}^{k}_l)-\frac{\theta_l^{k+1}}{2}\Vert W^{k+1}_l-\overline{W}_l^{k+1}\Vert^2_2 \ (\text{Equation \eqref{ineq:ineq2} and Equation \eqref{ineq:ineq4}})\\ \nonumber&=-(\nabla_{\overline{W}^{k+1}_l}\phi+\theta^{k+1}_l(W^{k+1}_l-\overline{W}^{k+1}_l))^T(W^{k}_l-W^{k+1}_l)-(\nabla_{\overline{W}^{k+1}_l}\phi)^T(W^{k+1}_l-{W}^{k}_l)-\frac{\theta_l^{k+1}}{2}\Vert W^{k+1}_l-\overline{W}_l^{k+1}\Vert^2_2 (\text{Lemma \ref{lemma:lemma 1}})\\\nonumber&=\frac{\theta_l^{k+1}}{2}\Vert W^{k+1}_l-\overline{W}_l^{k+1}\Vert^2_2+\theta^{k+1}_l(W^{k+1}_l-\overline{W}^{k+1}_l)^T(\overline{W}^{k+1}_l-{W}^{k}_l)\\&=\frac{\theta^{k+1}_l}{2}(\Vert W^{k+1}_l-W^k_l\Vert^2_2-\Vert \overline{W}^{k+1}_l-W^k_l\Vert^2_2)
\\&=\frac{\theta^{k+1}_l}{2}\Vert W^{k+1}_l-W^k_l\Vert^2_2 \ (\overline{W}^{k+1}_l=W^k_l).
\end{align*}
Let $\alpha^{k+1}_l=\theta^{k+1}_l$, then Equation \eqref{eq:mDLAM w optimality} still holds.
\end{proof}
\textbf{Proof of Lemma \ref{lemma:bounded objective}}
\begin{proof}
In Algorithm \ref{algo:mDLAM}:\\
(a). We sum Equation \eqref{eq:mDLAM w optimality},  Equation \eqref{eq:mDLAM z optimality} and Equation \eqref{eq:mDLAM a optimality} from $l=1$ to $L$ and from $k=0$ to $K$ to obtain
\begin{align}
    \nonumber&F(\textbf{W}^0,\textbf{z}^0,\textbf{a}^0)- F(\textbf{W}^K,\textbf{z}^K,\textbf{a}^K)\\&\geq\sum\nolimits_{k=0}^K(\sum\nolimits_{l=1}^L( \frac{\alpha_l^{k+1}}{2}\Vert W^{k+1}_l-W^k_l\Vert^2_2+ \frac{\gamma^{k+1}_l}{2}\Vert z^{k+1}_l-z^k_l\Vert^2_2)+\sum\nolimits_{l=1}^{L-1} \frac{\delta_l^{k+1}}{2}\Vert a^{k+1}_l-a^k_l\Vert^2_2).
    \label{ineq:mDLAM sequence convergence}
\end{align}
So $F(\textbf{W}^K,\textbf{z}^K,\textbf{a}^K)\leq F(\textbf{W}^0,\textbf{z}^0,\textbf{a}^0)$. This proves the upper boundness of $F$.
Let $K\rightarrow\infty$ in Equation \eqref{ineq:mDLAM sequence convergence}, since $F>0$ is lower bounded, we have
\begin{align}
&\sum\nolimits_{k=0}^K(\sum\nolimits_{l=1}^L( \frac{\alpha_l^{k+1}}{2}\Vert W^{k+1}_l-W^k_l\Vert^2_2+ \frac{\gamma^{k+1}_l}{2}\Vert z^{k+1}_l-z^k_l\Vert^2_2)+\sum\nolimits_{l=1}^{L-1} \frac{\delta_l^{k+1}}{2}\Vert a^{k+1}_l-a^k_l\Vert^2_2)<\infty.\label{eq:mDLAM finite sum}
\end{align}
Since the sum of this infinite series is finite, every term converges to 0.
This means that $\lim\nolimits_{k\rightarrow\infty}W^{k+1}_l-W^k_l=0$, $\lim\nolimits_{k\rightarrow\infty}z^{k+1}_l-z^k_l=0$ and $\lim\nolimits_{k\rightarrow\infty}a^{k+1}_l-a^k_l=0$. In other words, $\lim\nolimits_{k\rightarrow\infty}\textbf{W}^{k+1}-\textbf{W}^{k}=0$, $\lim\nolimits_{k\rightarrow\infty}\textbf{z}^{k+1}-\textbf{z}^{k}=0$, and $\lim\nolimits_{k\rightarrow\infty}\textbf{a}^{k+1}-\textbf{a}^{k}=0$.\\
(b). Because $F(\textbf{W}^k,\textbf{z}^k,\textbf{a}^k)$ is bounded, by the definition of coercivity and Assumption \ref{ass:coercive}, $(\textbf{W}^k,\textbf{z}^k,\textbf{a}^k)$ is bounded.\\ 
\end{proof}
\textbf{Proof of Lemma \ref{lemma:subgradient bound}}
\begin{proof}
As shown in Remark 2.2 in \cite{xu2013block}, 
\begin{align*}
    &\partial_{\textbf{W}^{k+1}} F=\{\partial_{W_1^{k+1}} F\}\times \{\partial _{W_2^{k+1}} F\}\times\cdots\times\{\partial _{W_L^{k+1}} F\},
\end{align*}
where $\times$ denotes Cartesian Product.\\
In Algorithm \ref{algo:mDLAM}, for $W^{k+1}_l$, according to Line 6 of Algorithm \ref{algo:mDLAM}, if\\ $F(\textbf{W}^{k+1}_{\leq l},\textbf{z}^{k+1}_{\leq l-1},\textbf{a}^{k+1}_{\leq l-1})< F(\textbf{W}^{k+1}_{\leq l-1},\textbf{z}^{k+1}_{\leq l-1},\textbf{a}^{k+1}_{\leq l-1})$, then \\
\begin{align}
    \partial_{W^{k+1}_l} \nonumber&F=\partial\Omega_l(W^{k+1}_l)+\nabla_{W_l^{k+1}} \phi(a^{k+1}_{l-1},W^{k+1}_l,z^{k+1}_l) \text{(Definition of $F$ in  Problem \ref{prob:inequality constrained DNN})}\\\nonumber&=\nabla_{W_l^{k+1}} \phi(a^{k+1}_{l-1},W^{k+1}_l,z^{k+1}_l)-\nabla_{\overline{W}^{k+1}_l} \phi(a^{k+1}_{l-1},\overline{W}^{k+1}_l,z^{k}_l)-\theta^{k+1}_l(W^{k+1}_l-\overline{W}^{k+1}_l)+\partial\Omega_l(W^{k+1}_l)\\\nonumber&+\nabla_{\overline{W}^{k+1}_l} \phi(a^{k+1}_{l-1},\overline{W}^{k+1}_l,z^{k}_l)+\theta^{k+1}_l(W^{k+1}_l-\overline{W}^{k+1}_l)\\\nonumber
    &=\rho(W^{k+1}_l-\overline{W}^{k+1}_l)a^{k+1}_{l-1}(a^{k+1}_{l-1})^T-\rho(z^{k+1}_l-z^k_l)(a^{k+1}_{l-1})^T-\theta^{k+1}_l(W^{k+1}_l-\overline{W}^{k+1}_l)+\partial\Omega_l(W^{k+1}_l)\\&+\nabla_{\overline{W}^{k+1}_l} \phi(a^{k+1}_{l-1},\overline{W}^{k+1}_l,z^{k}_l)+\theta^{k+1}_l(W^{k+1}_l-\overline{W}^{k+1}_l). \label{eq:w subgradient}
\end{align}
On one hand, we have
\begin{align}
\nonumber&\Vert\rho(W^{k+1}_l-\overline{W}^{k+1}_l)a^{k+1}_{l-1}(a^{k+1}_{l-1})^T-\rho(z^{k+1}_l-z^k_l)(a^{k+1}_{l-1})^T-\theta^{k+1}_l(W^{k+1}_l-\overline{W}^{k+1}_l)\Vert\\
\nonumber&\leq\rho\Vert(W^{k+1}_l-\overline{W}^{k+1}_l)a^{k+1}_{l-1}(a^{k+1}_{l-1})^T\Vert+\rho\Vert(z^{k+1}_l-z^k_l)(a^{k+1}_{l-1})^T\Vert+\theta^{k+1}_l\Vert W^{k+1}_l-\overline{W}^{k+1}_l\Vert \text{(Triangle Inequality)}\\\nonumber&\leq \rho\Vert W^{k+1}_l-\overline{W}^{k+1}_l\Vert\Vert a^{k+1}_{l-1}\Vert\Vert a^{k+1}_{l-1}\Vert+\rho\Vert z^{k+1}_l-z^k_l\Vert\Vert a^{k+1}_{l-1}\Vert+\theta^{k+1}_l\Vert W^{k+1}_l-\overline{W}^{k+1}_l\Vert \text{(Cauchy-Schwarz Inequality)}\\&
\leq \rho M_{\textbf{a}}\Vert z^{k+1}_l-z^{k}_l\Vert+(\rho M^2_{\textbf{a}}+\theta^{k+1}_{l})\Vert W^{k+1}_l-\overline{W}^{k+1}_l\Vert \ \text{(Lemma \ref{lemma:bounded objective})}\label{ineq:w subgradient bound}\\&\nonumber\leq\rho M_{\textbf{a}}\Vert z^{k+1}_l-z^{k}_l\Vert+(\rho M^2_{\textbf{a}}+\theta^{k+1}_{l})\Vert W^{k+1}_l-(W^{k}_l+\omega^k(W^{k}_l-W^{k-1}_l))\Vert(\text{Nesterov Acceleration})\\&\nonumber\leq\rho M_{\textbf{a}}\Vert z^{k+1}_l-z^{k}_l\Vert+(\rho M^2_{\textbf{a}}+\theta^{k+1}_{l})\Vert W^{k+1}_l-W^{k}_l\Vert+(\rho M^2_{\textbf{a}}+\theta^{k+1}_{l})\Vert W^{k}_l-W^{k-1}_l\Vert
(\text{Triangle Inequality and }\omega^k<1).
\end{align}
On the other hand,  the optimality condition of  Equation \eqref{eq:update W} yields
\begin{align*}
    0\in \partial\Omega_l(W^{k+1}_l)+\nabla_{\overline{W}^{k+1}_l} \phi(a^{k+1}_{l-1},\overline{W}^{k+1}_l,z^{k}_l)+\theta^{k+1}_l(W^{k+1}_l-\overline{W}^{k+1}_l).
\end{align*}
Therefore, there exists $g^{k+1}_{1,l}\in  \partial_{W^{k+1}_l} F$ such that 
\begin{align*}
    &\Vert g^{k+1}_{1,l}\Vert\leq \rho M_{\textbf{a}}\Vert z^{k+1}_l-z^{k}_l\Vert+(\rho M^2_{\textbf{a}}+\theta^{k+1}_{l})\Vert W^{k+1}_l-W^{k}_l\Vert+(\rho M^2_{\textbf{a}}+\theta^{k+1}_{l})\Vert W^{k}_l-W^{k-1}_l\Vert.
\end{align*}
This shows that there exists $g^{k+1}_1=g^{k+1}_{1,1}\times g^{k+1}_{1,2}\times\cdots\times g^{k+1}_{1,L}\in\partial_{\textbf{W}^{k+1}} F$ and $C_2=\max(\rho M_{\textbf{a}},\rho M^2_{\textbf{a}}+\theta^{k+1}_1,\rho M^2_{\textbf{a}}+\theta^{k+1}_2,\cdots,\rho M^2_{\textbf{a}}+\theta^{k+1}_L)$ such that
\begin{align}
    \Vert g^{k+1}_{l}\Vert&\leq C_2(\Vert \textbf{W}^{k+1}-\textbf{W}^{k}\Vert+\Vert \textbf{z}^{k+1}-\textbf{z}^{k}\Vert+\Vert \textbf{W}^{k}-\textbf{W}^{k-1}\Vert) \label{ineq:w bound 1}.
\end{align}
Otherwise, we have
\begin{align*}
&\Vert\rho(W^{k+1}_l-\overline{W}^{k+1}_l)a^{k+1}_{l-1}(a^{k+1}_{l-1})^T-\rho(z^{k+1}_l-z^k_l)(a^{k+1}_{l-1})^T-\theta^{k+1}_l(W^{k+1}_l-\overline{W}^{k+1}_l)\Vert\\&\leq\rho M_{\textbf{a}}\Vert z^{k+1}_l-z^{k}_l\Vert+(\rho M^2_{\textbf{a}}+\theta^{k+1}_{l})\Vert W^{k+1}_l-\overline{W}^{k+1}_l\Vert (\text{Equation \eqref{ineq:w subgradient bound}})\\&=\rho M_{\textbf{a}}\Vert z^{k+1}_l-z^{k}_l\Vert+(\rho M^2_{\textbf{a}}+\theta^{k+1}_{l})\Vert W^{k+1}_l-{W}^{k}_l\Vert (\overline{W}^{k+1}_l={W}^{k}_l).
\end{align*}
The optimality condition of  Equation \eqref{eq:update W} yields
\begin{align*}
    0\in \partial\Omega_l(W^{k+1}_l)+\nabla_{\overline{W}^{k+1}_l} \phi(a^{k+1}_{l-1},\overline{W}^{k+1}_l,z^{k}_l)+\theta^{k+1}_l(W^{k+1}_l-\overline{W}^{k+1}_l). 
\end{align*}
By Equation \eqref{eq:w subgradient}, we know that there exists $g^{k+1}_{1,l}\in  \partial_{W^{k+1}_l} F$ such that 
\begin{align}
    &\Vert g^{k+1}_{1,l}\Vert\leq \rho M_{\textbf{a}}\Vert z^{k+1}_l-z^{k}_l\Vert+(\rho M^2_{\textbf{a}}+\theta^{k+1}_{l})\Vert W^{k+1}_l-{W}^{k}_l\Vert. \label{ineq:w bound 2}
\end{align}
Combining Equation \eqref{ineq:w bound 1} with Equation \eqref{ineq:w bound 2}, we  show that there exists $g^{k+1}_1=g^{k+1}_{1,1}\times g^{k+1}_{1,2}\times\cdots\times g^{k+1}_{1,L}\in\partial_{\textbf{W}^{k+1}} F$ and $C_2=\max(\rho M_{\textbf{a}},\rho M^2_{\textbf{a}}+\theta^{k+1}_1,\rho M^2_{\textbf{a}}+\theta^{k+1}_2,\cdots,\rho M^2_{\textbf{a}}+\theta^{k+1}_L)$ such that
\begin{align*}
    \Vert g^{k+1}_{l}\Vert&\leq C_2(\Vert \textbf{W}^{k+1}-\textbf{W}^{k}\Vert+\Vert \textbf{z}^{k+1}-\textbf{z}^{k}\Vert+\Vert \textbf{W}^{k}-\textbf{W}^{k-1}\Vert).
\end{align*}
   \end{proof}
\textbf{Proof of Lemma \ref{lemma:objective decrease}}
\begin{proof}
We add Equation \eqref{eq:mDLAM w optimality},   Equation \eqref{eq:mDLAM z optimality}, and Equation \eqref{eq:mDLAM a optimality} from $l=1$ to $L$ to obtain
\begin{align*}
    &F(\textbf{W}^k,\textbf{z}^k,\textbf{a}^k)- F(\textbf{W}^{k+1},\textbf{z}^{k+1},\textbf{a}^{k+1})\\&\geq\sum\nolimits_{l=1}^L( \frac{\alpha_l^{k+1}}{2}\Vert W^{k+1}_l-W^k_l\Vert^2_2+ \frac{\gamma^{k+1}_l}{2}\Vert z^{k+1}_l-z^k_l\Vert^2_2)+\sum\nolimits_{l=1}^{L-1} \frac{\delta_l^{k+1}}{2}\Vert a^{k+1}_l-a^k_l\Vert^2_2.
\end{align*}
Let $C_5=\min(\frac{\alpha_l^{k+1}}{2},\frac{\gamma_l^{k+1}}{2},\frac{\delta_l^{k+1}}{2})>0$, we have
\begin{align}
  \nonumber &F(\textbf{W}^k,\textbf{z}^k,\textbf{a}^k)- F(\textbf{W}^{k+1},\textbf{z}^{k+1},\textbf{a}^{k+1})\\&\nonumber \geq C_5(\sum\nolimits_{l=1}^L(\Vert W^{k+1}_l-W^k_l\Vert^2_2+\Vert z^{k+1}_l-z^k_l\Vert^2_2)+\sum\nolimits_{l=1}^{L-1} \Vert a^{k+1}_l-a^k_l\Vert^2_2)\\&=C_5(\Vert \textbf{W}^{k+1}-\textbf{W}^{k}\Vert^2_2
  +\Vert \textbf{z}^{k+1}-\textbf{z}^{k}\Vert^2_2+\Vert \textbf{a}^{k+1}-\textbf{a}^{k}\Vert^2_2)\label{eq:objective decrease}
  \\\nonumber&\geq 0.
\end{align}
By Lemma \ref{lemma:bounded objective}(b) and a monotone sequence is convergent if it is bounded, then $F(\textbf{W}^k,\textbf{z}^k,\textbf{a}^k)$ is convergent.
\end{proof}
\textbf{Proof of Theorem \ref{thero:mDLAM convergence}}
 \begin{proof}
By Lemma \ref{lemma:bounded objective} (a), $\lim\nolimits_{k\rightarrow\infty}\textbf{W}^{k+1}-\textbf{W}^{k}=0$. By Lemma \ref{lemma:bounded objective} (b), there exists a subsequence $\textbf{W}^s$ such that $\textbf{W}^s\rightarrow \textbf{W}^*$, where $W^*$ is a limit point. From Lemma \ref{lemma:subgradient bound}, there exist $g_1^s\in \partial_{\textbf{W}^s} F$ such that $\Vert g_1^s\Vert \rightarrow 0$ as $s\rightarrow \infty$. According to the definition of limiting subdifferential, we have $0\in \partial{_{\textbf{W}^*}} F$. In other words, $\textbf{W}^*$ is a stationary point of $F$ in Problem \ref{prob:inequality constrained DNN}.
\end{proof}

   \textbf{Proof of Theorem \ref{thero:mDLAM convergence rate}}
\begin{proof}
In Algorithm \ref{algo:mDLAM}, we prove this by the KL Property.\\
Firstly, we consider Equation \eqref{eq:update z} and Equation \eqref{eq:update a}, by Lemma \ref{lemma:bounded objective},  $ h_l(\overline{z}^{k+1}_l-\nabla\phi_{\overline{z}^{k+1}_l}/\rho)-a^k_l$ and $ h_l(z^{k+1})-\overline{a}^{k+1}_l+\nabla_{\overline{a}^{k+1}_l}\phi/\tau^{k+1}_l$ are bounded, i.e. there exist constants $D_1$ and $D_2$ such that
\begin{align*}
    &\mid h_l(\overline{z}^{k+1}_l-\nabla_{\overline{z}^{k+1}_l}\phi/\rho)-a^k_l\mid< D_1.
    \\&\  \mid h_l(z^{k+1})-\overline{a}^{k+1}_l+\nabla_{\overline{a}^{k+1}_l}\phi/\tau^{k+1}_l\mid< D_2.
\end{align*}
Let $\varepsilon=\max(D_1,D_2)$, then the solutions to Equation \eqref{eq:update z} and Equation \eqref{eq:update a} are simplified as follows:
\begin{align}
    &z_l^{k+1}\leftarrow \overline{z}^{k+1}_{l}-\nabla_{\overline{z}^{k+1}_{l}}\phi/\rho.\label{eq:simplified z}\\&
    a_l^{k+1}\leftarrow \overline{a}^{k+1}_l-\nabla_{\overline{a}^{k+1}_l}\phi/\tau^{k+1}_l.\label{eq:simplified a}
\end{align}
This is because $h_l(z_l^{k+1})-\varepsilon\leq a^{k}_l\leq h_l(z_l^{k+1})+\varepsilon$ and $h_l(z_l^{k+1})-\varepsilon\leq a^{k+1}_l\leq h_l(z_l^{k+1})+\varepsilon$ hold in Equation \eqref{eq:update z} and Equation \eqref{eq:update a}, respectively. \\
Next, we prove that given $\varepsilon=\max(D_1,D_2)$, there exists $C_3=\max(\rho M^2_\textbf{W}+\tau_1^{k+1},\rho M^2_\textbf{W}+\tau_2^{k+1},\rho M^2_\textbf{W}+\tau_3^{k+1},\cdots,\rho M^2_\textbf{W}+\tau_{L-1}^{k+1}, 2\rho M_\textbf{W}M_\textbf{a}+\rho M_\textbf{z})$, some $ g^{k+1}_3\in \partial_{\textbf{z} ^{k+1}}F$ and $g^{k+1}_4\in \partial_{\textbf{a}^{k+1}} F$ such that
\begin{align*}
    &\Vert g^{k+1}_3\Vert=0,\\&
        \Vert g^{k+1}_4\Vert \leq C_3(\Vert \textbf{a}^{k+1}-\textbf{a}^{k}\Vert+\Vert\textbf{a}^{k}-\textbf{a}^{k-1}\Vert+\Vert\textbf{W}^{k+1}-\textbf{W}^{k}\Vert+\Vert \textbf{z}^{k+1}-\textbf{z}^{k} \Vert).
\end{align*}
As shown in \cite{wang2015global,xu2013block}, 
\begin{align*}
    &\partial_{\textbf{z}^{k+1}} F=\partial_{z_1^{k+1}} F\times \partial _{z_2^{k+1}} F\times\cdots\times\partial _{z_L^{k+1}} F,\\
    &\nabla _{\textbf{a}^{k+1}} F=\nabla_{a_1^{k+1}} F\times \nabla_{a_2^{k+1}} F\times\cdots\times \nabla_{a_{L-1}^{k+1}} F,\\
\end{align*}
where $\times$ denotes Cartesian Product.\\
For $z^{k+1}_l(l< L)$, according to Line 18 of Algorithm \ref{algo:mDLAM}, no matter \\$F(\textbf{W}^{k+1}_{\leq l},\textbf{z}^{k+1}_{\leq l},\textbf{a}^{k+1}_{\leq l-1})\geq F(\textbf{W}^{k+1}_{\leq l},\textbf{z}^{k+1}_{\leq l-1},\textbf{a}^{k+1}_{\leq l-1})$ or not, we have
\begin{align*}
    \partial_{z^{k+1}_l} F&=\nabla_{z^{k+1}_l} \phi(a^{k+1}_{l-1},W^{k+1}_l,z^{k+1}_l)\\&=\nabla_{z^{k+1}_l} \phi(a^{k+1}_{l-1},W^{k+1}_l,z^{k+1}_l)
    -\nabla_{\overline{z}^{k+1}_l} \phi(a^{k+1}_{l-1},W^{k+1}_l,\overline{z}^{k+1}_l)-\rho( z^{k+1}_l-\overline{z}^{k+1}_l) (\text{Equation \eqref{eq:simplified z}})\\&=0
.\end{align*}
For $z^{k+1}_L$, according to Line 12 of Algorithm \ref{algo:mDLAM}, no matter\\ $F(\textbf{W}^{k+1}_{\leq L},\textbf{z}^{k+1}_{\leq L},\textbf{a}^{k+1}_{\leq L-1})\geq F(\textbf{W}^{k+1}_{\leq L},\textbf{z}^{k+1}_{\leq L-1},\textbf{a}^{k+1}_{\leq L-1})$ or not, we have
\begin{align*}
    \partial_{z^{k+1}_L} F&=\nabla_{z^{k+1}_L} \phi(a^{k+1}_{L-1},W^{k+1}_L,z^{k+1}_L)+\partial R(z^{k+1}_L;y)\\&=\nabla_{z^{k+1}_L} \phi(a^{k+1}_{L-1},W^{k+1}_L,z^{k+1}_L)+\partial R(z^{k+1}_L;y)+\nabla_{\overline{z}^{k+1}_L} \phi(a^{k+1}_{L-1},W^{k+1}_L,\overline{z}^{k+1}_L)\\&+\rho( z_L-\overline{z}^{k+1}_L)-\nabla_{\overline{z}^{k+1}_L} \phi(a^{k+1}_{L-1},W^{k+1}_L,\overline{z}^{k+1}_L)-\rho( z^{k+1}_L-\overline{z}^{k+1}_L)\\&=\nabla_{z^{k+1}_L} \phi(a^{k+1}_{L-1},W^{k+1}_L,z^{k+1}_L)
    -\nabla_{\overline{z}^{k+1}_L} \phi(a^{k+1}_{L-1},W^{k+1}_L,\overline{z}^{k+1}_L)-\rho( z^{k+1}_L-\overline{z}^{k+1}_L)\\& (0\in\partial R(z^{k+1}_L;y)+\nabla_{\overline{z}^{k+1}_L} \phi(a^{k+1}_{L-1},W^{k+1}_L,\overline{z}^{k+1}_L)+\rho( z^{k+1}_L-\overline{z}^{k+1}_L)\text{by the optimality condition of Equation \eqref{eq:update zl}})\\&=0.
\end{align*}
Therefore, there exists $g^{k+1}_{3,l}= \nabla_{z^{k+1}_l} F$ such that $\Vert g^{k+1}_{3,l}\Vert=0$.
This shows that there exists $g^{k+1}_3=g^{k+1}_{3,1}\times g^{k+1}_{3,2}\times\cdots\times g^{k+1}_{3,L}=\nabla_{\textbf{z}^{k+1}}F$ such that 
\begin{align}
    \Vert g^{k+1}_3\Vert=0. \label{eq: z subgradient bound}
\end{align}
For $a^{k+1}_l$, we have
\begin{align*}
    \partial_{a^{k+1}_l} F&=\nabla_{a^{k+1}_l}\phi(a^{k+1}_l,W_{l+1}^k,z^{k+1}_{l+1})\\&=\nabla_{a^{k+1}_l}\phi(a^{k+1}_l,W^{k+1}_{l+1},z^{k+1}_{l+1})-\nabla_{\overline{a}^{k+1}_{l}} \phi(\overline{a}^{k+1}_l,W_{l+1}^k,z^{k}_{l+1})-\tau_l^{k+1}( a^{k+1}_{l}-\overline{a}^{k+1}_{l})(\text{Equation \eqref{eq:simplified a}})\\&=\rho (W^{k+1}_{l+1})^T(W^{k+1}_{l+1}a^{k+1}_l-z^{k+1}_{l+1})-\rho (W^{k}_{l+1})^T(W^{k}_{l+1}\overline{a}^{k+1}_l-z^{k}_{l+1})-\tau_l^{k+1}( a^{k+1}_{l}-\overline{a}^{k+1}_{l})\\&=\rho(W^{k+1}_{l+1})^T W^{k+1}_{l+1}(a^{k+1}_l-\overline{a}^{k+1}_l)+\rho(W^{k+1}_{l+1})^T (W^{k+1}_{l+1}-W^k_{l+1})\overline{a}^{k+1}_l\\&+\rho(W^{k+1}_{l+1}-W^{k}_{l+1})^T W^k_{l+1}\overline{a}^{k+1}_l-\rho(W^{k+1}_{l+1})^T(z^{k+1}_{l+1}-z^{k}_{l+1})-\rho(W^{k+1}_{l+1}-W^{k}_{l+1})^Tz^k_{l+1}-\tau_l^{k+1}( a^{k+1}_{l}-\overline{a}^{k+1}_{l}).
\end{align*}
Therefore
\begin{align*}
    \Vert\partial_{a^{k+1}_l} F\Vert &\leq \rho\Vert W^{k+1}_{l+1}\Vert \Vert W^{k+1}_{l+1}\Vert \Vert a^{k+1}_l-\overline{a}^{k+1}_l\Vert+\rho\Vert W^{k+1}_{l+1}\Vert \Vert W^{k+1}_{l+1}-W^k_{l+1}\Vert\Vert\overline{a}^{k+1}_l\Vert\\&+\rho\Vert W^{k+1}_{l+1}-W^{k}_{l+1}\Vert \Vert W^k_{l+1}\Vert\Vert\overline{a}^{k+1}_l\Vert+\rho\Vert W^{k+1}_{l+1}\Vert \Vert z^{k+1}_{l+1}-z^{k}_{l+1}\Vert\\&+\rho\Vert W^{k+1}_{l+1}-W^{k}_{l+1}\Vert \Vert z^k_{l+1}\Vert+\tau_l^{k+1}\Vert a^{k+1}_{l}-\overline{a}^{k+1}_{l}\Vert \\&(\text{Triangle Inequality and Cauthy-Schwarz Inequality})\\&\leq \rho M^2_\textbf{W}\Vert a^{k+1}_l-\overline{a}^{k+1}_l\Vert+\rho M_\textbf{W} \Vert W^{k+1}_{l+1}-W^k_{l+1}\Vert M_\textbf{a}+\rho\Vert W^{k+1}_{l+1}-W^{k}_{l+1}\Vert M_\textbf{W}M_\textbf{a}+\rho M_\textbf{W} \Vert z^{k+1}_{l+1}-z^{k}_{l+1}\Vert\\&+\rho\Vert W^{k+1}_{l+1}-W^{k}_{l+1}\Vert M_\textbf{z}+\tau_l^{k+1}\Vert a^{k+1}_{l}-\overline{a}^{k+1}_{l}\Vert \ \text{(Lemma \ref{lemma:bounded objective})}\\&=(\rho M^2_\textbf{W}+\tau_l^{k+1})\Vert a^{k+1}_l-\overline{a}^{k+1}_l\Vert+(2\rho M_\textbf{W}M_\textbf{a}+\rho M_\textbf{z}) \Vert W^{k+1}_{l+1}-W^k_{l+1}\Vert+\rho M_\textbf{W} \Vert z^{k+1}_{l+1}-z^{k}_{l+1}\Vert.
\end{align*}
According to Line 22 of Algorithm \ref{algo:mDLAM}, if \\$ F(\textbf{W}^{k+1}_{\leq l},\textbf{z}^{k+1}_{\leq l},\textbf{a}^{k+1}_{\leq l})< F(\textbf{W}^{k+1}_{\leq l},\textbf{z}^{k+1}_{\leq l},\textbf{a}^{k+1}_{\leq l-1})$, then we have
\begin{align*}
    \Vert\partial_{a^{k+1}_l} F\Vert&\leq (\rho M^2_\textbf{W}+\tau_l^{k+1})\Vert a^{k+1}_l-{a}^{k}_l-({a}^{k}_l-{a}^{k-1}_l)\omega^k\Vert+(2\rho M_\textbf{W}M_\textbf{a}+\rho M_\textbf{z}) \Vert W^{k+1}_{l+1}-W^k_{l+1}\Vert+\rho M_\textbf{W} \Vert z^{k+1}_{l+1}-z^{k}_{l+1}\Vert \\& (\text{Nestrov Acceleration})\\&\leq(\rho M^2_\textbf{W}+\tau_l^{k+1})\Vert a^{k+1}_l-{a}^{k}_l\Vert+(\rho M^2_\textbf{W}+\tau_l^{k+1})\Vert a^k_l-a^{k-1}_l\Vert+(2\rho M_\textbf{W}M_\textbf{a}+\rho M_\textbf{z}) \Vert W^{k+1}_{l+1}-W^k_{l+1}\Vert\\&+\rho M_\textbf{W} \Vert z^{k+1}_{l+1}-z^{k}_{l+1}\Vert \ (\text{Triangle Inequality and $\omega^k<1$}).
\end{align*}
Therefore, there exists $g^{k+1}_{4,l}\in \partial_{a^{k+1}_l} F$ such that
\begin{align}
    \nonumber\Vert g^{k+1}_{4,l}\Vert&\leq (\rho M^2_\textbf{W}+\tau_l^{k+1})\Vert a^{k+1}_l-{a}^{k}_l\Vert+(\rho M^2_\textbf{W}+\tau_l^{k+1})\Vert a^k_l-a^{k-1}_l\Vert+(2\rho M_\textbf{W}M_\textbf{a}+\rho M_\textbf{z}) \Vert W^{k+1}_{l+1}-W^k_{l+1}\Vert\\&+\rho M_\textbf{W} \Vert z^{k+1}_{l+1}-z^{k}_{l+1}\Vert. \label{eq:a bound1}
\end{align}
Otherwise, 
\begin{align*}
    \Vert \partial_{a^{k+1}_l} F\Vert&\leq(\rho M^2_\textbf{W}+\tau_l^{k+1})\Vert a^{k+1}_l-{a}^{k}_l\Vert+(2\rho M_\textbf{W}M_\textbf{a}+\rho M_\textbf{z}) \Vert W^{k+1}_{l+1}-W^k_{l+1}\Vert+\rho M_\textbf{W} \Vert z^{k+1}_{l+1}-z^{k}_{l+1}\Vert
     \ (\overline{a}^{k+1}_l={a}^{k}_l). 
\end{align*}
Therefore, there exists $g^{k+1}_{4,l}\in \partial_{a^{k+1}_l} F$ such that
\begin{align}
     \Vert g^{k+1}_{4,l} \Vert&\leq(\rho M^2_\textbf{W}+\tau_l^{k+1})\Vert a^{k+1}_l-{a}^{k}_l\Vert+(2\rho M_\textbf{W}M_\textbf{a}+\rho M_\textbf{z}) \Vert W^{k+1}_{l+1}-W^k_{l+1}\Vert+\rho M_\textbf{W} \Vert z^{k+1}_{l+1}-z^{k}_{l+1}\Vert. \label{eq:a bound2}
\end{align}
Combining Equation \eqref{eq:a bound1} and Equation \eqref{eq:a bound2}, we  show that there exists $g^{k+1}_4=g^{k+1}_{4,1}\times g^{k+1}_{4,2}\times\cdots\times g^{k+1}_{4,L}\in\partial_{\textbf{a}^{k+1}} F$ and $C_3=\max(\rho M^2_\textbf{W}+\tau_1^{k+1},\rho M^2_\textbf{W}+\tau_2^{k+1},\rho M^2_\textbf{W}+\tau_3^{k+1},\cdots,\rho M^2_\textbf{W}+\tau_{L-1}^{k+1}, 2\rho M_\textbf{W}M_\textbf{a}+\rho M_\textbf{z})$ such that
\begin{align}
    \Vert g^{k+1}_4\Vert \leq C_3(\Vert \textbf{a}^{k+1}-\textbf{a}^{k}\Vert+\Vert\textbf{a}^{k}-\textbf{a}^{k-1}\Vert+\Vert\textbf{W}^{k+1}-\textbf{W}^{k}\Vert+\Vert \textbf{z}^{k+1}-\textbf{z}^{k} \Vert). \label{eq: a subgradient bound}
\end{align}
Combining Lemma \ref{lemma:subgradient bound}, Equation \eqref{eq: z subgradient bound} and Equation \eqref{eq: a subgradient bound}, we prove that there exists $g^{k+1}\in \partial F(\textbf{W}^{k+1},\textbf{z}^{k+1},\textbf{a}^{k+1})=\{ \partial_{\textbf{W}^{k+1}}F,\partial_{\textbf{z}^{k+1}}F,\partial_{\textbf{a}^{k+1}}F\}$ and $C_4=\max(C_2,C_3,\rho)$ such that
\begin{align}
    \Vert g^{k+1}\Vert&\leq C_4(\Vert \textbf{a}^{k+1}-\textbf{a}^{k}\Vert+\Vert\textbf{a}^{k}-\textbf{a}^{k-1}\Vert+\Vert\textbf{W}^{k+1}-\textbf{W}^{k}\Vert+\Vert\textbf{W}^{k}-\textbf{W}^{k-1}\Vert+\Vert \textbf{z}^{k+1}-\textbf{z}^{k} \Vert).
    \label{eq:subgradient bound}
\end{align}
Finally, we prove the linear convergence rate by the KL Property given Equation \eqref{eq:subgradient bound} and Equation \eqref{eq:objective decrease}.
Because $F$ is locally strongly convex with a constant $\mu$, $F$ satisfies the KL Property by Lemma \ref{lemma:locally strongly convex KL property}. Let $F^*=F(\textbf{W}^*,\textbf{z}^*,\textbf{a}^*)$ be the convergent value of $F$, by Lemma \ref{lemma:objective decrease},  $F(\textbf{W}^k,\textbf{z}^k,\textbf{a}^k)\rightarrow F^*$, then for any $\eta_1>0$ there exists $k_2\in \mathbb{N}$ such that it holds for $k>k_2$ that $     F^*<F(\textbf{W}^k,\textbf{z}^k,\textbf{a}^k)<F^*+\eta_1$.
Also by Lemma \ref{lemma:bounded objective}(a) and  Equation \eqref{eq:subgradient bound}, $g^{k+1}\rightarrow 0$ as $k\rightarrow \infty$, then for any $\eta_2>0$ there exists $k_3\in \mathbb{N}$, such that it holds for $k>k_3$ that $\Vert g^{k+1}\Vert<\eta_2$. Therefore, for any $k>k_1=\max(k_2,k_3)$, $(\textbf{W}^k,\textbf{z}^k,\textbf{a}^k)\in\{(\textbf{W},\textbf{z},\textbf{a}):|F^*<F(\textbf{W},\textbf{z},\textbf{a})<F^*+\eta_1\cap \exists g\in F(\textbf{W},\textbf{z},\textbf{a}) \ s.t. \ \Vert g\Vert<\eta_2\}$. By the KL Property and Lemma \ref{lemma:locally strongly convex KL property}, it holds that
\begin{align*}
    1&\leq \Vert g^{k+1}\Vert/(\mu \sqrt{F(\textbf{W}^{k+1},\textbf{z}^{k+1},\textbf{a}^{k+1})-F^*}) \\&\leq C_4(\Vert \textbf{a}^{k+1}-\textbf{a}^{k}\Vert+\Vert\textbf{a}^{k}-\textbf{a}^{k-1}\Vert+\Vert\textbf{W}^{k+1}-\textbf{W}^{k}\Vert+\Vert\textbf{W}^{k}-\textbf{W}^{k-1}\Vert+\Vert \textbf{z}^{k+1}-\textbf{z}^{k} \Vert)/(\mu \sqrt{F(\textbf{W}^{k+1},\textbf{z}^{k+1},\textbf{a}^{k+1})-F^*}) \\&  (\text{Equation \eqref{eq:subgradient bound}})\\&\leq C^2_4(\Vert \textbf{a}^{k+1}-\textbf{a}^{k}\Vert+\Vert\textbf{a}^{k}-\textbf{a}^{k-1}\Vert+\Vert\textbf{W}^{k+1}-\textbf{W}^{k}\Vert+\Vert\textbf{W}^{k}-\textbf{W}^{k-1}\Vert+\Vert \textbf{z}^{k+1}-\textbf{z}^{k} \Vert)^2/(\mu^2 (F(\textbf{W}^{k+1},\textbf{z}^{k+1},\textbf{a}^{k+1})-F^*))\\&\leq (5C^2_4(\Vert \textbf{a}^{k+1}-\textbf{a}^{k}\Vert^2_2+\Vert\textbf{a}^{k}-\textbf{a}^{k-1}\Vert^2_2+\Vert\textbf{W}^{k+1}-\textbf{W}^{k}\Vert^2_2+\Vert\textbf{W}^{k}-\textbf{W}^{k-1}\Vert^2_2+\Vert \textbf{z}^{k+1}-\textbf{z}^{k} \Vert^2_2))/(\mu^2 (F(\textbf{W}^{k+1},\textbf{z}^{k+1},\textbf{a}^{k+1})-F^*)) \\& (\text{Mean Inequality})\\&\leq (5C^2_4(F(\textbf{W}^{k-1},\textbf{z}^{k-1},\textbf{a}^{k-1})-F(\textbf{W}^{k+1},\textbf{z}^{k+1},\textbf{a}^{k+1})))/(C_5\mu^2 (F(\textbf{W}^{k+1},\textbf{z}^{k+1},\textbf{a}^{k+1})-F^*)) (\text{Equation \eqref{eq:objective decrease}}).
\end{align*}
This indicates that
\begin{align*}
    &(C_5\mu^2+5C^2_4)(F(\textbf{W}^{k+1},\textbf{z}^{k+1},\textbf{a}^{k+1})-F^*)\leq 5C^2_4 (F(\textbf{W}^{k-1},\textbf{z}^{k-1},\textbf{a}^{k-1})-F^*).
\end{align*}
Let $0<C_1=\frac{5C^2_4}{C_5\mu^2+5C^2_4}<1$, we have
\begin{align*}
    F(\textbf{W}^{k+1},\textbf{z}^{k+1},\textbf{a}^{k+1})-F^*\leq C_1(F(\textbf{W}^{k-1},\textbf{z}^{k-1},\textbf{a}^{k-1})-F^*).
\end{align*}
    So in summary,  for any $\rho$, there exist $\varepsilon=\max(D_1,D_2)$, $k_1=\max(k_2,k_3)$, and $0<C_1=\frac{5C^2_4}{C_5\mu^2+5C^2_4}<1$ such that
   \begin{align*}
    F(\textbf{W}^{k+1},\textbf{z}^{k+1},\textbf{a}^{k+1})-F^*\leq C_1(F(\textbf{W}^{k-1},\textbf{z}^{k-1},\textbf{a}^{k-1})-F^*).
\end{align*}
    for $k>k_1$. In other words, the linear convergence rate is proven.
\end{proof}

\end{document}